\documentclass[review,hidelinks,onefignum,onetabnum]{siamart251104}
\usepackage[T1]{fontenc} 

\usepackage{lipsum}
\usepackage{amsfonts}
\usepackage{graphicx}
\usepackage{epstopdf}
\usepackage{algorithmic}
\ifpdf
  \DeclareGraphicsExtensions{.eps,.pdf,.png,.jpg}
\else
  \DeclareGraphicsExtensions{.eps}
\fi

\usepackage{amssymb}
\usepackage{amsmath}
\allowdisplaybreaks
\usepackage{cancel}
\usepackage{graphicx} 
\usepackage{animate} 
\usepackage{subcaption}
\usepackage[export]{adjustbox}
\usepackage{xcolor}

\usepackage{hyperref} 

\newsiamremark{remark}{Remark}
\newsiamremark{hypothesis}{Hypothesis}
\newsiamremark{example}{Example}
\newsiamremark{examples}{Examples}
\crefname{hypothesis}{Hypothesis}{Hypotheses}
\newsiamthm{claim}{Claim}
\newsiamthm{assumption}{Assumption}
\newsiamremark{fact}{Fact}
\crefname{fact}{Fact}{Facts}

\headers{Initial Value Learning Based MFGs on Graphs}{Yaxin Feng, Yang Xiang, and Haomin Zhou}

\title{Discrete Mean Field Games on Finite Graphs As Initial Value Optimization}

\author{Yaxin Feng\footnotemark[1]\thanks{Department of Mathematics, The Hong Kong University of Science and Technology,
Hong Kong SAR, China (\email{yfengba@connect.ust.hk, maxiang@ust.hk}).}
\and Yang Xiang$^*$\footnotemark[2]\thanks{HKUST Shenzhen-Hong Kong Collaborative Innovation Research Institute, Futian, Shenzhen, China.}
\and Haomin Zhou\thanks{School of Mathematics, Georgia Institute of Technology, GA 30332 USA (\email{hmzhou@math.gatech.edu}).}
}

\usepackage{amsopn}

\usepackage[margin=1in]{geometry}

\begin{document}\nolinenumbers

\maketitle

\begin{abstract}
In this paper, we propose an initial value fomulation of the discrete mean field games on finite graphs (Graph MFG), and design a neural network based approach to solve it. Graph MFG describes infinite, non-cooperative and interactive homogeneous agents move on node states through the edges to optimize their own goals. Nash Equilibrium of the Graph MFG is characterized by a coupled ordinary differential equations (ODE) system, including the discrete forward continuity equation and the discrete backward Hamilton-Jacobi equation. 
In this paper, we mainly focus on the potential mean field games (Potential MFG) on finite graphs, which has an infinite-dimensional constrained optimization structure. We reformulate Potential MFG as an initial value finite-dimentional optimization problem with dynamics constrains, names Graph MFG-IV. Specifically, the initial condition of the Hamilton-Jacobi equation is regarded as the unique variable, constrained by the coupled Hamilton-Jacobi and continuity equation system as the ODE integrator. This formulation is a reduced-order model, which avoids time-discretization of the infinite-dimensional path and has a much smaller searching space than the general path-wise problem setting. We design a neural network-based approach to solve the Graph MFG-IV problem. 
\end{abstract}



\section{Introduction}\label{sec:introduction}



Mean field games (MFG) describe an infinite number of non-cooperative agents optimizing their own goals within an interacting crowd~\cite{lasry2007mean,huang2006large}.
It can be regarded as a game between an indistinguishable representative agent and a generic agent (population).
It is successfully applied to many areas, including game theory~\cite{cardaliaguet2021introduction,chew2016potential}, economics~\cite{achdou2022income,gueant2010mean,carmona2020applications}, opinion dynamics~\cite{liu2022deep,bauso2016opinion}, robotics~\cite{liu2018mean, elamvazhuthi2019mean, kang2020joint}, epidemic modeling~\cite{aurell2022optimal,roy2023recent} and so on. It have been found to be a powerful tool in generative modeling, being connected to a closely related topic optimal transport (OT)~\cite{ruthotto2020machine,zhang2023mean} and Schrödinger Bridge~\cite{liu2022deep}.  Recently, neural networks have been proposed for efficiency and reduction to solve MFG~\cite{hu2023recent,lawal2022physics,lin2021alternating,huang2023bridging,han2024learning,carmona2022convergence}.
The variational structure of MFG under the dynamics constraint of continuity equation is:
\begin{equation}\label{eq:Potential MFG-opt}
\begin{aligned}
\operatorname*{inf}_{\rho, v} \quad & J(\rho, v) = \int_0^T \mathcal{K}(\rho,v) + \mathcal{F}(\rho(t)) \,dt + \mathcal{G}(\rho(T)) \\
\text{subject to} \quad & \partial_t \rho + \nabla \cdot (\rho v) = 0, \quad \rho(0, x) = \mu_0(x),
\end{aligned}
\end{equation}
Here, $\mathcal{K}$ is the transportation energy \textit{w.r.t.} the velocity field $v:\mathbb{R}^d\rightarrow\mathbb{R}^d$ and the population density $\rho \in \mathcal{P}(\mathbb{R}^d)$, in which $\mathcal{P}(\mathbb{R}^d)$ is the space of all probability densities. Besides, $\mathcal{F},\mathcal{G}$ are the potentials \textit{w.r.t.} $\rho$.  
This structure can be named as potential mean field games (Potential MFG). At Nash equilibrium, no agent can decrease its cost by modifying the strategy (control, \textit{e.g.} velocity)
unilaterally, the actual density $\rho \in \mathcal{P}(\mathbb{R}^d)$ satisfies the continuity equation or Fokker Planck equation, and the minimum cost of the single agent is the value function, which satisfies the Hamilton-Jacobi (HJ) equation. $\mathcal{F}$ is the weighted sum of the potentials \textit{w.r.t.} population $\rho$, and
the terminal energies $\mathcal{G}$ to be the density estimation functionals, such as point-wise $L_1$ distance
or the Kullback-Leibler (KL) divergence. 
When $\mathcal{F}=0$ and the terminal condition strictly satisfies $\rho(T)=\mu_T$, Potential MFG becomes the dynamical OT defined by Benamou-Brenier~\cite{benamou2000computational}, also known 
as the Wasserstein-2 metric $\mathcal{W}_2^2(\mu_0, \mu_T)$. 


Many mass transport problems can be described as Wasserstein Hamiltonian flows (WHF)~\cite{chow2020Wasserstein}. The WHF is a coupled density and value dynamics system, 
including optimal transport, Schrödinger Bridge, Schrödinger equations, mean field models, etc. The first equation is generally a Kolmogorov equations, \textit{e.g.} continuity equation or Fokker Planck equation, and the latter one is a HJ equation. MFG can also be described as WHF, with the initial condition of density function and the terminal condition of value function. 

In this paper, we design a computational strategy to solve the discrete Potential MFG on finite graphs, in which the population dynamics and the value functions evolve on the discrete nodes that connected by undirected edges. 

Discrete mean field games on the finite graphs (named as Graph MFG) 
has important applications in many areas, such as social networks, 
opinion dynamics, traffic and routing problems~\cite{yang2017learning,calderone2017markov,huang2021dynamic,cabannes2021solving,shou2022multi,chen2023learning}, etc. The solution existence and uniqueness of Graph MFG are discussed in~\cite{gueant2011infinity,gueant2015existence}.
In these studies, the dynamics of the agents is mainly described in Markov decision process, so the mass transport equations are the Markov chains, and the velocity field is not explicitly defined. This problem is conventionally solved by dynamics programming~\cite{huang2021dynamic,cabannes2021solving} and fixed-point iteration~\cite{tanaka2020linearly}. Recently, reinforcement learning algorithms, such as actor-critic methods, have also been designed to solve these Markov decision process~\cite{chen2023learning,yang2017learning}.


There are some related studies on the discrete mass transport (\textit{e.g.} gradient flow, optimal transport, etc.) on finite graphs or states~\cite{chow2012fokker,maas2011gradient,mielke2011gradient}. \cite{chow2012fokker,chow2017entropy} inllustrate that the Fokker Planck equation on graphs can be described as the gradient flow of free potential functional or the Markov process. 
Benamou-Brenier formulation of OT on graphs is proposed to solve the discrete Schrödinger equation~\cite{chow2019discrete}. \cite{gangbo2024well} analyzes the well-poseness of the initial value HJ equation on graphs, which can describe the optimal control problems, and \cite{cui2025finite} focuses on numerical analysis of two proposed finite difference schemes of the initial value problem of HJ equation on graphs. 

Graph MFG is generally defined as a path-wise, infinite-dimensional optimization problem. In this paper, we propose a reformulation of the Potential MFG on graphs as an initial value optimization problem from the viewpoint of WHF, and design a neural networks-based algorithm for it.

This paper is organized as follows. In
\cref{sec:GMFG}, the Graph MFG is defined. A reformulation to the initial value Graph MFG is proposed in \cref{sec:GMFG-inv}. We design a neural networks-based algorithm, and illustrate a warm start scheme for the network training in \cref{sec:alg}. The experimental results are presented in \cref{sec:experiments}, and a conclusion follows in
\cref{sec:conclusions}.

\section{Discrete Mean Field Games on Finite Graphs}
\label{sec:GMFG}

We consider an undirected, connected graphs $G=(E,V,\mathbf{W})$, with no self-loops or multiple edges~\cite{gross2003handbook,chow2019discrete}. We define $V=\{a_i\}_{i=1}^n$ to be the vertex set, $(i,j) \in E$ as the edge between nodes $a_i$ and $a_j$ on the edge set $E$, $N(i)=\left\{a_j \in V:(i, j) \in E\right\}$ as the neighbor set of node $a_i$, $i=1,\cdots,n$, and $w_{ij} \in \mathbf{W}$ as the weight of the edge $(i, j)\in E$ satisfying $w_{ij}=w_{ji} > 0$.

\subsection{Notations}

Here, we define some notations on graphs, following~\cite{chow2017entropy}.

The probability set (simplex) supported on the finite graph $G$ is defined as 
$$\mathcal{P}(G)=\{\left(\rho_i\right)_{i=1}^n \in \mathbb{R}^n \mid \sum_{i=1}^n \rho_i=1, \rho_i \geq 0, \quad \text { for any } i \in V\},$$
where $\rho_i$ is the probability on the node $i$. The interior of $\mathcal{P}(G)$ is denoted by $\mathcal{P}_o(G)$, in which $\rho_i > 0$ on each node.

The vector field $v:V\times V \rightarrow \mathbb{R}$ represents the velocity on the graph $G$, which is a skew-symmetric matrix defined on the edge set $E$: 
$$v_{ij}=\begin{cases}
  -v_{ji} \quad \text{ if } (v_{ij})_{(i,j) \in E}\\0 \qquad \text{ otherwise}
&\end{cases}.$$

The discrete flux $m=\rho v,$ and $m:V\times V \rightarrow \mathbb{R}$ on the finite graph $G$ is an anti-symmetric function such that $m_{ij}=-m_{ji},$ $$m:=\left(v_{ij} \theta_{ij}(\rho)\right)_{(i, j) \in E},$$ where $\theta_{ij}(\rho)$ being the weighted function of $\rho_i$ and $\rho_j$.
In this work, the weight function $\theta_{ij}(\cdot)$ of density $\rho$ is defined as $\theta_{ij}(\rho) = \frac{\rho_i+\rho_j}{2}$. Noted that there are also other types of weight functions listed in Remark~\ref{remk:theta_rho}.

The discrete value function $S=(S_i)_{i=1}^n :V\rightarrow \mathbb{R}$ induces a vector field $$(\nabla_\text{G}S)_{ij}:=\sqrt{w_{ij}}(S_i-S_j)_{(i,j) \in E},$$ where $w_{ij}$'s are the weights on the edges, and $\nabla_\text{G}S:V\times V \rightarrow \mathbb{R}.$

Given two vector fields $v$ and $\tilde{v}$, the discrete inner product of the vector fields on the finite graph $G$ and $\rho \in \mathcal{P}(G)$ is $$\langle v,\tilde{v}\rangle_\rho=\frac{1}{2}\sum_{(i,j) \in E}v_{ij}\tilde{v}_{ij}\theta_{ij}(\rho),$$ in which the factor $\frac{1}{2}$ is needed because all the values defined on the edges are counted twice. We have $\langle v,v\rangle_\rho=\frac{1}{2}\sum_{i,j \in E}v_{ij}^2\theta_{ij}(\rho)$ when $\tilde{v}=v.$
We choose the transport cost to be kinetic energy $\mathcal{K}(\rho,v)=\frac{1}{2} \langle v, v \rangle_\rho$.

Divergence of $\rho v$ on G is
$$\text{div}_\text{G}(\rho v) := -\left(\sum_{j \in N(i)} \sqrt{w_{ij}}v_{ij} \theta_{ij}(\rho)\right)_{i=1}^n.$$

An integration by part formula on graphs is
\begin{equation}\label{eq:dint_by_part}
\begin{aligned}
-\sum_{i=1}^n \text{div}_\text{G} (\rho v)|_i \xi_i &= \sum_{i=1}^n \sum_{j \in N(i)}  \sqrt{w_{ij}}v_{ij} \theta_{ij}(\rho) \xi_i \\
&= \frac{1}{2}(\sum_{(i,j) \in E}  \sqrt{w_{ij}}v_{ij} \xi_i \theta_{ij}(\rho) + \sum_{(j,i) \in E}  \sqrt{w_{ji}}v_{ji} \xi_j \theta_{ji}(\rho)) \\
&= \frac{1}{2}\sum_{(i,j) \in E}  \sqrt{w_{ij}}v_{ij} \cdot (\xi_i-\xi_j) \theta_{ij}(\rho) \\
&= \langle\nabla_\text{G} \xi, v\rangle_\rho,
\end{aligned}
\end{equation}
and $\sum_{i=1}^n \text{div}_\text{G} (\rho v) |_i = 0$ when letting $\xi_i=1, i\in V.$

\begin{remark}\label{remk:theta_rho}
Noted that the choice of the average weighted density $\theta_{ij}(\mathbf\rho)=\theta_{ij}^A(\mathbf\rho)=\frac{\rho_i+\rho_j}{2}$ is not unique. There are other types of weighted functions, such as the upwind weight $\theta_{ij}^U(\mathbf\rho)=\rho_i, S_i < S_j$, the logarithmic weight $\theta_{ij}^L(\mathbf\rho)=\frac{\rho_i-\rho_j}{\log \left(\rho_i\right)-\log \left(\rho_j\right)}$~\cite{maas2011gradient}, and the arithmetic mean $\theta^{AM}_{ij}(\rho)=\frac{1}{2}\left(\frac{\rho_i}{d_i}+\frac{\rho_j}{d_j}\right)$, with $d_i=\frac{\sum_{j \in N(i)} w_{i j}}{\sum_{i=1}^n \sum_{j \in N(i)} w_{i j}}$ representing the volume at node $i$~\cite{chow2022dynamical,cui2022time}. 
\end{remark}



\subsection{Discrete Mean Field Games on Finite Graphs}


\begin{definition}[Potential Mean Field Games on Finite Graphs]
\label{def:pw-Potential MFG-d}
Define the Graph MFG with the energy $\mathcal{K}: \left\{(\rho, v) \mid \rho \in \mathcal{P}({G}), v \in L^2(V\times V;\theta(\rho)) \right\} \to \mathbb{R}$, the potentials $\mathcal{F},\mathcal{G}:\mathcal{P}(G)\rightarrow \mathbb{R}$ on the graph $G=(E,V,\mathbf{W})$ as:
\begin{equation}\label{eq:GMFG-opt}
\begin{aligned}
\operatorname*{inf}_{\rho, v} \quad & J_1(\rho, v) = \int_0^T 
\mathcal{K}(\rho(t),v(t)) + \mathcal{F}(\rho(t)) \,dt +  \mathcal{G}(\rho(T))
\end{aligned}
\end{equation}
subject to
\begin{equation}\label{eq:ce}
\begin{aligned}
&\frac{d \rho_i}{dt}-\sum_{j \in N(i)}\sqrt{w_{ij}}v_{ij}\theta_{ij}(\rho)=0, \\
& \rho_i(0) = {\mu_0}_i, \quad \text{ for } i=1,\cdots,n, j\in N(i).
\end{aligned}
\end{equation}
where $\rho_i(t)$ and the $v_{ij}(t)$ are the probability function and the velocity field at time $t$ on node $i$ and edge $(i,j), i,j=1,\cdots,n,$ respectively.
\end{definition}

The constrained optimization in \cref{def:pw-Potential MFG-d} is an infinite dimensional path-wise problem, for which we name it Graph MFG-PW for convenience. 

In this paper, 
we consider potential $\mathcal{F}:=\lambda_V\mathcal{V} + \lambda_W\mathcal{W} + \lambda_U\mathcal{U} 
$. Specifically, the discrete linear potential $\mathcal{V}(\rho):=\sum_{i=1}^n \mathbb{V}_i \rho_i$, the discrete interaction potential of the agents $\mathcal{W}(\rho):=\frac{1}{2} \sum_{j=1}^n \sum_{i=1}^n \mathbb{W}_{ij} \rho_i \rho_j$, in which $\mathbb{W}_{ij}$ are symmetric and positive definite matrices, 
and the discrete nonlinear potential $\mathcal{U}(\rho)$ is chosen to be the discrete negative Boltzmann Shannon entropy $\mathcal{B}(\rho):=\sum_{i=1}^n \rho_i \log\rho_i,$ weighted by $\lambda_B$. The discrete terminal energy $\mathcal{G}(\rho(T))$ in this work is based on the $L_1$ distance 
$\mathcal{G}_{L1}({\rho(T)}):=\|{\rho(T)}-{\mu_T}\|_1$
or the KL divergence $\mathcal{G}_{\mathrm{KL}}({\rho(T)}):=\sum_{i=1}^n{\rho_i(T)}\log \frac{{\rho_i(T)}}{{\mu_T}_i}$, weighted by $\lambda_G$.
These potenetials are all convex \textit{w.r.t.} $\rho$ and $\rho(T)$. We remark that there are other discrete potentials, such as the discrete Fisher information $\mathcal{I}(\rho):=\frac{1}{2} \sum_{(i, j) \in E} \omega_{ij}\left(\log \rho_i-\log \rho_j\right)^2 \theta_{ij}(\rho)$ that can also be included.

Define $f({\rho}):=\lambda_V V({\rho})+\lambda_W W({\rho})
+\lambda_U U({\rho})$ to be the cost of the agent, and $g({\rho(T)})$ to be the terminal energy of it, $f,g:\mathcal{P}(G)\rightarrow \mathbb{R}$. 
When these cost terms are the first variation of the potentials in~\cref{eq:GMFG-opt}, $$
\begin{aligned}
&f(\rho)=\frac{\delta \mathcal{F}(\rho)}{\delta \rho} = \lambda_V\frac{\delta \mathcal{V}(\rho)}{\delta \rho} + \lambda_W\frac{\delta \mathcal{W}(\rho)}{\delta \rho} + \lambda_U\frac{\delta \mathcal{U}(\rho)}{\delta \rho}, \\
&g(\rho(T))=\frac{\delta \mathcal{G}(\rho(T))}{\delta \rho(T)},
\end{aligned}
$$ in which 
$
\left. V_i(\rho)=\frac{\delta \mathcal{V}(\rho)}{\delta \rho} \right|_i=\mathbb{V}_i,  \left.W_i(\rho)=\frac{\delta \mathcal{W}(\rho)}{\delta \rho}\right|_i = \sum_{j=1}^n \mathbb{W}_{ij}\rho_j$, and $U(\rho)=\frac{\delta \mathcal{U}(\rho)}{\delta \rho}$ is chosen as $B_i(\rho)=\left.\frac{\delta \mathcal{B}(\rho)}{\delta \rho}\right|_i=\log\rho_i + 1, i=1,\cdots,n,\lambda_U=\lambda_B.$ Then we can obtain the necessary first-order optimality condition as well as the Nash equilibrium of~\cref{eq:GMFG-opt}. 


\begin{definition}[Nash Equilibrium of the Discrete Mean Field Games on Finite Graphs]
\label{def:hj-ce-d}
Nash equilibrium of the Graph MFG in \cref{def:pw-Potential MFG-d} satisfies the discrete continuity equation and the discrete HJ equation system on the graph $G$:
\begin{equation}\label{eq:bw-hj-ce}
  \begin{cases}
    \frac{d \rho_i}{dt} -\sum_{j \in N(i)}\sqrt{w_{ij}}v_{ij} \theta_{ij}(\rho) = 0, \\
    -\frac{d S_i}{dt} + \frac{1}{2} \sum_{j \in {N(i)}} w_{ij}(S_i-S_j)^2\frac{\partial{\theta_{ij}(\rho)}}{\partial \rho_i} = f_i({\rho}), \\ \rho_i(0) = {\mu_0}_i, \quad S_i(T) = {g}_i({\rho(T)}), \quad \text{for } i=1,\cdots,n, \\ 
    v_{ij} = \sqrt{w_{ij}}(S_j-S_i), \quad j\in N(i),
  \end{cases}
\end{equation}
where $\rho_i(t)$ and the $S_i(t)$ are the probability function and the value function at time $t$ on node $i$, the $i=1,\cdots,n$. 
\end{definition}


Denote the augmented Lagrangian \begin{equation}\label{eq:L1}
\mathcal{L}_1(\rho,v)=J_1(\rho,v)-\int_0^T \sum_{i=1}^n \Phi_i[\frac{d \rho_i}{dt} -\sum_{j \in N(i)}\sqrt{w_{ij}}v_{ij} \theta_{ij}(\rho)]\,dt, 
\end{equation}
where $\Phi_i$ is the Lagrangian multiplier. The Nash equilibrium in~\cref{eq:bw-hj-ce} is the first-order necessary optimality condition of Graph MFG in~\cref{def:pw-Potential MFG-d}, which can be proved by $\delta \mathcal{L}_1(\rho,v)=0.$

Let $\delta \mathcal{F}(\rho)=\sum_{i=1}^n f_i(\rho)\delta\rho_i$ and $\delta^2 \mathcal{F}(\rho)=\sum_{i=1}^n f'_i(\rho) (\delta\rho_i)^2$ denote the first and second variations of $\mathcal{F}(\rho)$ for convenience. We use the same representation for $\mathcal{G}(\rho(T))$. 

\begin{assumption}\label{ass:convexity}
  Assume $\mathcal{F}(\rho)$ and $\mathcal{G}({\rho(T)})$ are convex functionals in the Graph MFG. 
  The second variations $\delta^2 \mathcal{F}(\rho)>0$ and $\delta^2 \mathcal{G}({\rho(T)})>0$.
\end{assumption}

\begin{remark}\label{remk:Lipschitz}

  Given the initial condition, the~\cref{eq:fw-hj-ce} system has a unique solution within $t \in [0,T)$ once it is locally Lipschitz continuous, and solution always exists when uniformly Lipschitz continuous.
  The weighted density function is $\theta_{ij}(\rho)=\frac{\rho_i+\rho_j}{2}$ here, thus when $\mathcal{F} = 0,$ ~\cref{eq:fw-hj-ce} is uniformly Lipschitz continuous. When $\mathcal{F} \neq 0,$ $\frac{\delta \mathcal{V}}{\delta \rho}$ and $\frac{\delta \mathcal{W}}{\delta \rho}$ with bounded interaction 
  kernel $\sup_{(i,j)\in E}\mathbb{W}_{ij}$ are uniformly Lipschitz continuous, and $\frac{\delta \mathcal{B}}{\delta \rho}$
  is locally Lipschitz continuous. 
\end{remark}

\begin{remark}\label{remk:exist-unique}
  The solution existence and uniqueness of the continuous and discrete mean field games are extensively discussed on~\cite{lasry2007mean,gomes2010discrete,gomes2013continuous,gueant2011infinity,gueant2015existence}. Specifically, mean field games on graphs have been discussed on~\cite{gueant2011infinity,gueant2015existence}.
\end{remark}

\section{Initial Value Optimization: Reformulation of Potential Mean Field Games on Finite Graphs}
\label{sec:GMFG-inv}

We reformulate the infinite-dimensional path-wise Graph MFG-PW to a finite-dimensional initial value optimization defined below. 

\begin{definition}[Potential Mean Field Games on Finite Graphs as Initial Value Optimization]
\label{def:iv-Potential MFG-d}
Define the initial value optimization of Graph MFG on $G$ (Graph MFG-IV):
\begin{equation}\label{eq:GMFG-iv} 
\begin{aligned}
\operatorname*{min}_{S_0} \quad & J_2(S_0) = \int_0^T \sum_{i=1}^n\frac{1}{4}\sum_{j\in N(i)} w_{ij}(S_i-S_j)^2 \theta_{ij}(\rho) + \mathcal{F}(\rho(t)) \,dt + \mathcal{G}(\rho(T)) 
\end{aligned}
\end{equation}
subject to
\begin{equation}\label{eq:fw-hj-ce}
\begin{cases}
\frac{d \rho_i}{dt}-\sum_{j \in N(i)}w_{ij}(S_j-S_i)\theta_{ij}(\rho)=0, \\
-\frac{dS_i}{dt}+\frac{1}{2}\sum_{j\in N(i)}w_{ij}(S_i-S_j)^2 \frac{\partial\theta_{ij}(\rho)}{\partial\rho_i}=
{f}_i(\rho)\\
\rho_i(0) = {\mu_0}_i(x), \quad S_i(0) = {S_0}_i, \quad \text{for } i=1,\cdots,n
\end{cases}
\end{equation}
where $\rho_i(t)$ and $S_i(t)$ are the probability function and the value function at time $t$ on node $i$ and edge $(i,j)$, $i,j=1,\cdots,n,$ respectively. 

\end{definition}


In the \cref{eq:GMFG-iv}, the objective function $J_2(S_0)$ has almost the same expression to the $J_1(\rho,v)$ in \cref{eq:GMFG-opt} except the variables: $S_0$ with finite dimension $n$ is the unique varibale in $J_2(S_0)$, while 
$\rho$ and $v$ in $J_1(\rho,v)$ are two time-dependent variables that are infinite-dimensional. 
This formulation of the initial value problem has a much smaller searching space than that of the path-wise formulation.

Denote the cost of the Grpah MFG-IV as $\tilde{J}_2(\rho,S):=J_2(S(0)),$ in which the path $(\rho,S)$ starts from the initial point $S(0)$.
Accordingly, the augmented Lagrangian about the generated path is 
\begin{equation}
  \begin{aligned}
  \tilde{\mathcal{L}}_2(\rho,S)=&\tilde{J}_2(\rho,S)-\int_0^T\sum_{i=1}^n{\phi}_i[\frac{d \rho_i}{dt}-\sum_{j \in N(i)}\sqrt{w_{ij}}(S_j-S_i)\theta_{ij}(\rho)]\,dt \\
  &\quad -\int_0^T\sum_{i=1}^n{\psi}_i[\frac{dS_i}{dt}-\frac{1}{2}\sum_{j\in N(i)}w_{ij}(S_i-S_j)^2 \frac{\partial\theta_{ij}(\rho)}{\partial\rho_i}+
  {f}_i(\rho)]\,dt,
  \end{aligned}
\end{equation}
where ${\phi}$ and ${\psi}$ are Lagrangian multipliers.

\begin{lemma}\label{le:hit.terminal}
  Assume that $f'(\rho)$ is Lipschitz continuous. If the path $(\rho^*,S^*)$ starts from the minimizer of Graph MFG-IV $S_0^*$ and evolves according to~\cref{eq:fw-hj-ce}, then it arrives the terminal condition of the HJ equation $S^*(T)=g(\rho^*(T))$. 
\end{lemma}
Proof of~\cref{le:hit.terminal}: 

Assume the solution of Graph MFG-IV is $S_0^*$, and the path generated from~\cref{eq:fw-hj-ce} is $(\rho^*,S^*)$. Assume a small perturbed $\delta S_0$ on the initial condition is of the order $O(\varepsilon)$, where $\varepsilon$ is a small value. Consider the path perturbed from the optimal initial point $S_0^*+\delta S_0$. The generated perturbed path following~\cref{eq:fw-hj-ce} is $(\rho^*+\delta \rho,S^*+\delta S),$ $\delta \rho$ and $\delta S$ are also of the order $O(\varepsilon)$, because $T$ is finite.

The perturbation satisfies the coupled linearized system from~\cref{eq:fw-hj-ce}:
\begin{equation}\label{eq:delta.HJ-CE}
\begin{aligned}  \left\{ \begin{array}{lr}
  \frac{d\delta \rho_i}{dt} - {\displaystyle\sum_{j \in N(i)}}{w_{ij}}(2(S_j^*-S_i^*)\frac{\partial \theta_{ij}(\rho^*)}{\partial\rho_i}\delta \rho_i + \theta_{ij}(\rho^*) (\delta S_j-\delta S_i))=0, 
  (\delta \text{CE})
  & \\ \frac{d\delta S_i}{dt}-\displaystyle\sum_{j \in N(i)}w_{ij}\frac{\partial \theta_{ij}(\rho^*)}{\partial \rho_i}(S_i^* - S_j^*) \cdot (\delta S_i-\delta S_j) + 
  f_i'(\rho^*)\delta \rho_i
  & \\-\frac{1}{2}(S_i^*-S_j^*)^2(\frac{\partial^2 \theta_{ij}(\rho^*)}{\partial \rho_i^2}+\frac{\partial^2 \theta_{ij}(\rho^*)}{\partial \rho_i \partial \rho_j}){\delta \rho_i} =0,  (\delta \text{HJE})
\end{array} \right. \end{aligned}
\end{equation} 
where $\delta S(0)=\delta S_0, \delta \rho(0)=0$. Here, $\delta \text{CE}$ and $\delta \text{HJ}$ are the first order perturbation of continuity equation and HJ equation in~\cref{eq:fw-hj-ce}, respectively. 

The first-order perturbation of $\tilde{\mathcal{L}_2}$ gives that 
\begin{equation}\label{eq:delta_L_2}
  \begin{aligned}
\delta \tilde{\mathcal{L}}_2 &= \delta \tilde{J}_2 - \langle \phi,\delta \text{CE} \rangle - \langle \psi,\delta \text{HJE} \rangle \\
&=\int_0^T \sum_{(i,j)\in E}[\frac{d\phi_i}{dt}+f_i(\rho^*)+w_{ij}(\phi_i-\phi_j)(S_j^*-S_i^*)\frac{\partial\theta_{ij}(\rho^*)}{\partial \rho_i}+\frac{1}{2}\psi_if_i'(\rho^*)\\
&\quad +\frac{1}{2}w_{ij}\|S_i^*-S_j^*\|^2\frac{\partial \theta_{ij}(\rho^*)}{\partial \rho_i}+\frac{1}{2}w_{ij}\psi_i(S_i^*-S_j^*)^2(\frac{\partial^2\theta_{ij}(\rho^*)}{\partial\rho_i\partial\rho_j}\\
&\quad+\frac{\partial^2\theta_{ij}(\rho^*)}{\partial\rho_i^2})]\delta \rho_i \,dt +\textstyle\sum_{i=1}^n[g(\rho_i^*(T)-\phi_i(\rho(T)))]\delta \rho_i(T)\\
&\quad +\int_0^T \sum_{(i,j)\in E}[\frac{d\psi_i}{dt}+w_{ij}(S_i^*-S_j^*)\theta_{ij}(\rho^*)-w_{ij}(\phi_i-\phi_j)\theta_{ij}(\rho^*)\\
&\quad +2w_{ij}\psi_i(S_i^*-S_j^*)\frac{\partial \theta_{ij}(\rho^*)}{\partial{\rho}_i}]\delta S_i \,dt +\textstyle\sum_{i=1}^n[\psi_i(0)\delta S_i(0)-\psi_i(T)\delta S_i(T)].
  \end{aligned}
\end{equation}
Here, $\phi(\rho(0))\delta\rho(0)=0$ because $\delta\rho(0)=0$.

The first-order optimality condition of $\tilde{\mathcal{L}_2}$ gives:
\begin{equation}\label{eq:kkt.delta_CE-HJ}
  \begin{aligned}  \left\{ \begin{array}{lr}
    \left.\frac{\delta\tilde{\mathcal{L}}_2}{\delta \rho}\right|_i=\frac{d\phi_i}{dt}+\displaystyle\sum_{j\in N(i)}\textstyle w_{ij}(\frac{1}{2}(S_i^*-S_j^*)-\phi_i+\phi_j)(S_i^*-S_j^*) \textstyle {\frac{\partial \theta_{ij}(\rho^*)}{\partial \rho_i}}
    \\\qquad +f_i(\rho^*)
    +\frac{1}{2}w_{ij}\psi_i(S_i^*-S_j^*)^2(\textstyle {\frac{\partial^2 \theta_{ij}(\rho^*)}{\partial \rho_i^2}}+\textstyle {\frac{\partial^2 \theta_{ij}(\rho^*)}{\partial \rho_i \partial \rho_j}})+\frac{1}{2}\psi_i f_i'(\rho^*)=0,
    & \\ \left.\frac{\delta\tilde{\mathcal{L}}_2}{\delta S}\right|_i=\frac{d\psi_i}{dt}+\displaystyle\sum_{j\in N(i)}w_{ij}(S_i^*-S_j^*)\theta_{ij}(\rho^*)-w_{ij}(\phi_i-\phi_j)\theta_{ij}(\rho^*)
    \\\qquad +2w_{ij}\psi_i(S_i^*-S_j^*)\frac{\partial \theta_{ij}(\rho^*)}{\partial \rho_i}=0,
  \end{array} \right. \end{aligned}
  \end{equation} 
and the boundary conditions 
\begin{equation}\label{eq:kkt.boundaries}
  \frac{\delta\tilde{\mathcal{L}}_2}{\delta \rho(T)}=g(\rho^*(T))-\phi(T)=0,\,\frac{\delta\tilde{\mathcal{L}}_2}{\delta S(0)}=\psi(0)=0,\,\frac{\delta\tilde{\mathcal{L}}_2}{\delta S(T)}=\psi(T)=0.
\end{equation}
When we let $\phi=S^*,\psi=0$, $\phi_i(t)$ and $\psi_i(t)$ satisfy the two ODEs in the system~\cref{eq:kkt.delta_CE-HJ}, with the initial condition $\phi(0)=S_0^*,\psi(0)=0$. Since~\cref{eq:kkt.delta_CE-HJ} is a Lipschitz continuous ODE system up to time $T$, the solution of the initial value problem of~\cref{eq:kkt.delta_CE-HJ} is unique.
The remaining boundary conditions in~\cref{eq:kkt.boundaries} $\phi(T)=S^*(T)=g(\rho^*(T)),\psi(T)=0$ are also satisfied. 
Therefore, the first-order optimality condition of $\tilde{\mathcal{L}}_2$ can be satisfied when $\phi=S^*,\psi=0$, and the terminal condition $S^*(T)=g(\rho^*(T))$ is satisfied.

Now, we show the relationship between our proposed finite-dimensional Graph MFG-IV and the infinite-dimensional Graph MFG-PW.
\begin{theorem}
\label{thm:reformulation}
If the solution of Graph MFG-IV (in \cref{def:iv-Potential MFG-d}) exists up to time $T$, the trajectory generated by~\cref{eq:fw-hj-ce} is a minimizer of the formulation Graph MFG-PW (in \cref{def:pw-Potential MFG-d}). Conversely, if the minimizer of Graph MFG-PW exists, its value function $S$ at $t=0$ provides a minimizer of Graph MFG-IV.

\end{theorem}
Proof of~\cref{thm:reformulation}.

\textit{(I) From Graph MFG-IV to Graph MFG-PW.}


In~\cref{le:hit.terminal}, we prove that the path of Graph MFG-IV arrives at the same terminal condition as in~\cref{eq:bw-hj-ce}.

Define the path of the solution of Graph MFG-IV to be $(\rho^*,S^*)$ and it gives the path $(\rho^*,v^*)$ that starts from $v(0)=v_0^*$, which is determined by $S(0)=S_0^*$. We want to show that $(\rho^*,v^*)$ is the solution of Graph MFG-PW.

Assume $\delta \rho$ and $\delta v$ are of the order $O(\varepsilon)$, $\varepsilon$ is a small value. Consider the augmented Lagrangian $\mathcal{L}_1(\rho,v)$ in~\cref{eq:L1} associated the path $(\rho^* + \delta \rho,v^* + \delta v)$:
\begin{equation}\label{eq:delta L_1}
\begin{aligned}
&\mathcal{L}_1(\rho^* + \delta \rho, v^* + \delta v) \\
=& 
\mathcal{K}(\rho^*+\delta \rho, v^*+\delta v)
+ \int_0^{T}\mathcal{F}(\rho^*+\delta \rho) \,dt + \mathcal{G}((\rho^*+ \delta \rho)(T)) - T_1,
\end{aligned}
\end{equation}
where 
$$
T_1 
=\int_0^{T}\sum_{i=1}^n \Phi_i[\frac{d{\rho_i^*}}{dt} + \frac{d\delta {\rho_i}}{dt}- \sum_{j \in N(i)}\sqrt{w_{ij}}(\theta_{ij}(\rho^*)+2\frac{\partial \theta_{ij}(\rho^*)}{\partial \rho_i}\delta \rho_i)(v_{ij}^*+\delta v_{ij})] \,dt,
$$
in which $\Phi_i$ is the Lagrangian multiplier, 
$$
\begin{aligned}
&\mathcal{K}(\rho^*+\delta \rho, v^*+\delta v) \\
=& \int_0^{T}\sum_{(i,j)\in E} \frac{1}{4}(\theta_{ij}(\rho^*)+2\frac{\partial \theta_{ij}(\rho^*)}{\partial \rho_i}\delta \rho_i)(\|v_{ij}^*\|^2 + 2v_{ij}^*\cdot \delta v_{ij} + \|\delta v_{ij}\|^2) \,dt.
\end{aligned}
$$

Organizing the terms according to the orders of $\varepsilon$ in \cref{eq:delta L_1}:
\begin{equation}\label{eq:delta L_1-1}
\begin{aligned}
&\mathcal{L}_1(\rho^*+\delta \rho,v^*+\delta v) \\
=&\mathcal{L}_1(\rho^*,v^*) +\{\textstyle{\int_0^{T}}\sum_{(i,j)\in E}[\textstyle{\frac{1}{2}}(v_{ij}^*\cdot \delta v_{ij})\textstyle{\theta_{ij}(\rho^*)} + \frac{1}{2}\frac{\partial \theta_{ij}(\rho^*)}{\partial \rho_i}\delta \rho_i \|v_{ij}^*\|^2 \\
&+f_i(\rho^*){\delta \rho_i}] \,dt + \textstyle\sum_{i=1}^n g_i({\rho^*(T)}){\delta {\rho_i(T)}} \\
& - \textstyle\sum_{i=1}^n \Phi_i [\textstyle{\frac{d\delta {\rho_i}}{dt}} - \sum_{j \in N(i)}\sqrt{w_{ij}}(2v_{ij}^*\frac{\partial \theta_{ij}(\rho^*)}{\partial\rho_i}\delta \rho_i + \theta_{ij}(\rho^*) \delta v_{ij})] \}  \\
& + \{ \textstyle{\int_0^{T}}\sum_{(i,j)\in E} [\textstyle{\frac{1}{4}}\|\delta v_{ij}\|^2 \theta_{ij}(\rho^*) + (v_{ij}^* \cdot \delta v_{ij})\textstyle{\frac{\partial \theta_{ij}(\rho^*)}{\partial\rho_i}}\delta \rho_i \\
& + 2\Phi_i \sqrt{w_{ij}}(\textstyle{\frac{\partial \theta_{ij}(\rho^*)}{\partial \rho_i}}\delta \rho_i \delta v_{ij}) +\frac{1}{2}f_i'(\rho^*){(\delta \rho_i)^2}] \,dt \\
& +\textstyle\sum_{i=1}^n\frac{1}{2} g_i'(\rho^*(T)){(\delta {\rho_i(T)})^2} \} 
+ \mathcal{O}(\varepsilon^3),
\end{aligned}
\end{equation}
the first and second order variation terms $\delta \mathcal{L}_1, \frac{1}{2}\delta^2 \mathcal{L}_1$ are seperated in two brackets here.

According to the first-order term in~\cref{eq:delta L_1-1} that can be calculated by integration by part:\begin{equation}\label{eq:delta_L1}
\begin{aligned}
  \delta \mathcal{L}_1 
  &=\int_0^{T}\frac{1}{2}\sum_{(i,j)\in E} \theta_{ij}(\rho^*)(v_{ij}^*+\sqrt{w_{ij}}(\Phi_i-\Phi_j))\cdot \delta v_{ij}\\
  &\quad +\sum_{i=1}^n[(\frac{d\Phi_i}{dt} + \sum_{j \in N(i)}(\frac{1}{2}\|v_{ij}^*\|^2 +\sqrt{w_{ij}}(v_{ij}^*(\Phi_i-\Phi_j)))\frac{\partial \theta_{ij}(\rho^*)}{\partial\rho_i}) \\
  &\quad + f_i(\rho^*)]\delta \rho_i \,dt + \sum_{i=1}^n (g_i(\rho^*(T))-{\Phi_i(T)})\delta {\rho_i(T)}.
\end{aligned}
\end{equation}
Here we have used ${\Phi_i(0)}\delta \rho_i(0)=0.$

If we choose the Lagrangian multiplier to be $\Phi=S^*,$ and let 
$v_{ij}^* = \sqrt{w_{ij}}(S_j^* - S_i^*)$, we have $\frac{\delta \mathcal{L}_1}{\delta v}=0$ in the first line. 
In the second line, according to the HJ equation in~\cref{eq:fw-hj-ce} of Graph MFG-IV that $S^*$ satisfies, we have $\frac{\delta \mathcal{L}_1}{\delta \rho}=0$. 
In this case, $\Phi^*(T)=S^*(T)=g(\rho^*(T))$ that induces $\frac{\delta \mathcal{L}_1}{\delta \rho(T)}=0$ in the third line.  
These mean that the path $(\rho^*,S^*)$ starting from the solution of Graph MFG-IV, $S_0^*$, satisfying~\cref{eq:fw-hj-ce}, 
makes $\delta \mathcal{L}_1(\rho^*,v^*)=\delta J_1(\rho^*,v^*)=0$, \textit{i.e.} $(\rho^*,v^*)$ is the critical point of~\cref{eq:GMFG-opt} in Graph MFG-PW. 
Next, we consider the second-order variation of the augmented Lagrangian $\mathcal{L}_1$ with $\delta \rho \neq 0, \delta v \neq 0$, \textit{i.e.} the third term in \cref{eq:delta L_1-1}, with $\Phi^*=S^*$, $v_{ij}^*=\sqrt{w_{ij}}(S_j^*-S_i^*)$:
\begin{equation}\label{eq:delta2 L_1}
\begin{aligned}
\frac{1}{2}\delta^2 \mathcal{L}_1 
&=\int_0^{T}\sum_{i=1}^n[ \frac{1}{4}\|\delta v_{ij}\|^2 \theta_{ij}(\rho^*) + \sum_{j\in N(i)}(v_{ij}^* \cdot \delta v_{ij})\frac{\partial \theta_{ij}(\rho^*)}{\partial\rho_i}\delta \rho_i \\
& \quad + \frac{1}{2}f_i'(\rho^*){(\delta \rho_i)^2}]  \,dt + T_2 + \sum_{i=1}^n\frac{1}{2}g_i'({\rho^*(T)}){(\delta \rho_i(T))^2},
\end{aligned}
\end{equation}
where
$$
\begin{aligned}
T_2 &= \int_0^{T}2S_i^* \sum_{(i,j)\in E}\sqrt{w_{ij}}(\frac{\partial \theta_{ij}(\rho^*)}{\partial \rho_i}\delta \rho_i \delta v_{ij}) \,dt\\
&=\int_0^{T}\sum_{(i,j)\in E} \sqrt{w_{ij}}(S_i^*-S_j^*) (\frac{\partial \theta_{ij}(\rho^*)}{\partial \rho_i}\delta \rho_i \delta v_{ij}) \,dt \\
&=-\int_0^{T}\sum_{(i,j)\in E}(v_{ij}^* \cdot \delta v_{ij})\frac{\partial \theta_{ij}(\rho^*)}{\partial\rho_i}\delta \rho_i \,dt,
\end{aligned}
$$
in which the second equation is from~\cref{eq:dint_by_part}, while $v_{ij}$, $\rho_i$ and $\xi_i$ are replaced by $\delta v_{ij}$, $\frac{\partial \theta_{ij}(\rho^*)}{\partial \rho_i}\delta \rho_i$ and $S_i^*,$ and the third one is according to the first-order condition that $v_{ij}^* = \sqrt{w_{ij}}(S_j^* - S_i^*)$.
Therefore, the second and third terms in the $\frac{1}{2}\delta^2\mathcal{L}_1$~\cref{eq:delta2 L_1} are canceled out:
\begin{equation}
\begin{aligned}
\delta^2 \mathcal{L}_1
&=\langle \delta v,\delta v \rangle_\rho + \int_0^{T} \sum_{i=1}^n f_i'(\rho^*){(\delta \rho_i)^2} \,dt + \sum_{i=1}^n g_i'(\rho^*(T)){(\delta \rho_i(T))^2}. 
\end{aligned}
\end{equation}
Since $\langle \delta v,\delta v \rangle_\rho>0$ and $f_i'(\rho^*)>0,$
$g_i'(\rho^*(T))>0$ because of the convexity of $\mathcal{F}(\rho)$ and $\mathcal{G}(\rho(T))$ from Assumption~\ref{ass:convexity}, we have $\delta^2 \mathcal{L}_2 > 0.$ 

Combining the critical point condition $\delta \mathcal{L}_1=0$ and $\delta^2 \mathcal{L}_1 > 0$, the path $(\rho^*,S^*)$ starting from the soluiton of Graph MFG-IV $S_0^*$ gives the solution of Graph MFG-PW $(\rho^*,v^*)$. 


\textit{(II) From Graph MFG-PW to Graph MFG-IV.} 

Assume the solution of Graph MFG-PW (in~\cref{def:pw-Potential MFG-d}) is $(\hat{\rho},\hat{v})$, \textit{i.e.} $J_1(\hat\rho,\hat{v}) = \inf_{\rho,v}J_1(\rho,v)$. We need to prove that when $\hat v=-\nabla \hat S$, $S(0)=\hat{S}(0)$ is the minimizer of $J_2$ in Graph MFG-IV (in~\cref{def:iv-Potential MFG-d}).

The solution of Graph MFG-PW $(\hat{\rho},\hat{v})$ satisfies the necessary first order optimality condition~\cref{eq:bw-hj-ce}, which is 
the same ODE system with the dynamics constraints~\cref{eq:fw-hj-ce} of Graph MFG-IV.
In other words, $(\hat \rho,\hat S)$ satisfies all the constraints in Graph MFG-IV. 
Therefore, when $\hat{S}(0)$ is a feasible initial point for Graph MFG-IV, the generated path is in the feasible set of Graph MFG-PW.
Since $\inf_{\rho,v} J_1(\rho, v)\leq \inf_{S_0} J_2(S_0)$ is always true, $J_2(\hat S(0))$ must be the minimum of Graph MFG-IV.

\section{Methodology}
\label{sec:alg}

To solve the Graph MFG, we design a neural network approach based on our new formulation Graph MFG-IV presented in~\cref{sec:GMFG-inv}. 


\textbf{Main Algorithm.} We use neural network to parameterize the initial condition of HJ equation $S_0 = {S_0}_{NN}(\Theta)$, where $\Theta$ is the set of parameters of the neural network. 
The path $(\rho,S)$ is calculated by an ODE integrator for the coupled HJ and continuity equation in~\cref{eq:fw-hj-ce} starting from ${S_0}_{NN}(\Theta),$ with an integration step $\Delta t$. The neural network is trained using a loss function based on a discretized form of~\cref{eq:GMFG-iv} with the step $\Delta t$. 
After the training is converged, we get the optimized $S_0^*$ and the solution of Graph MFG, $(\rho^*,S^*).$ The algorithm details are given in Alg.~\ref{alg:gmfgiv}.

\begin{algorithm}
\caption{(Graph MFG-IV) Neural Network for Discrete Mean Field games on Finite Graphs: Initial Value Optimization}
\label{alg:gmfgiv}
\begin{algorithmic}
\STATE{Given $\mu_0, \text{ guess } S_0 = {S_0}_{NN}(\Theta)
$}
\FOR{$k=0,\cdots,K$ epochs}
\FOR{$m=0,\cdots,M-1$ steps$, \forall i \in V$}
\STATE{Update $S$: $S_i^{(m+1)} := S_i^{(m)} + \frac{1}{2} \sum_{j \in {N(i)}} (w_{ij}(S_i^{(m)}-S_j^{(m)})^2\frac{\partial{\theta_{ij}(\rho)}}{\partial \rho_i^{(m)}} - \left.\frac{\delta\mathcal{F}(\rho^{(m)})}{\delta \rho}\right|_i)\,\Delta t$}
\STATE{Update $\rho$: $\rho_i^{(m+1)} := \rho_i^{(m)} + (\sum_{j \in N(i)}w_{ij}(S_j^{(m)}-S_i^{(m)}) \theta_{ij}(\rho^{(m)}))\,\Delta t$ }
\STATE{Collect the path during the ODE integration 
$(\rho^{(m+1)},S^{(m+1)}).$
}
\ENDFOR

Calculate loss function according to~\cref{eq:GMFG-iv}: $$J_2^{[k]}({S_0}_{NN}) = \sum_{m=0}^M [\frac{1}{2} \langle v, v \rangle_\rho + \mathcal{F}(\rho(t))]\Delta t + \mathcal{G}(\rho(T))$$

Update ${S_0}_{NN}(V;\Theta^{[k+1]})$ by gradient descent: $$\Theta^{[k+1]} := \Theta^{[k]} - \eta \nabla_{\Theta}J_2^{[k]}({S_0}_{NN})$$
\ENDFOR

\RETURN $S_0^* = {S_0}_{NN}(\Theta^{[K]})$
\end{algorithmic}
\end{algorithm}






The neural network for the Graph MFG-IV method has the following advantages:
(1) \textit{Reduce Order.} Here, shadow layer neural networks are enough to give good results. The number of network parameters is much smaller than that of nodes on a graph.
(2) \textit{Auto-Differentiation.} Chain rule gradients along the whole ODE integration \textit{w.r.t.} the initial condition $S_0$ are conveniently calculated by auto-differentiation in the neural network.
(3) \textit{Easier Training and Better Results.} Based on numerical results, training neural network to obtain the parameterized ${S_0}_{NN}$ is easier than directly optimizing the ${S_0}_i, i=1,\cdots,n$ on each node without using neural network, and can achieve better results.

\textbf{Warm-Start Scheme.}
A warm-start scheme is designed to further speed up and stablize the training, in which the initial parameters are trained based on an intuitive velocity field according to the initial density $\mu_0$ and the target density $\mu_T$. We represent the intuitive velocity field as the gradient of the intuitive potential. 
There are many intuitive velocity fields can be chosen, such as that straightly pointing from the center of the given initial density $\mu_0$ to the target density $\mu_T$, named as $\nabla P_1$; and the gradient of the difference between two densities, \textit{i.e.} $\nabla P_2=\nabla (\mu_0-\mu_T)$. We apply a linear combination of these two guessed velocity fields, with a parameter $\alpha\in[0,1]$. The gradient of the target intuitive potential is 
\begin{equation}\label{eq:warmup}
  \nabla P_{target}=\alpha \nabla P_1 + (1-\alpha)\nabla P_2.
\end{equation}

We train the initial set of parameters in ${S_0}_{InitialNet}$, such that $\nabla{S_0}_{InitialNet}$ is an approximation of $\nabla{P}_{target}$. We use mean square error $L_{mse}$ for this training. After convergence, we set the parameters of ${S_0}_{InitialNet}$ to be the initial values of the parameters in ${S_0}_{NN}$ in the main training. This warm-start algorithm is shown in Alg.~\ref{alg:warm-start}.


\begin{algorithm}
  \caption{Warm Start of the Initial Neural Network ${S_0}_{InitialNet}(V;\Theta_{init})$}
  \label{alg:warm-start}
  \begin{algorithmic}
    \STATE{Guess an intuitive potential $P_{target}$}
    \FOR{$\tilde{k}=0,\cdots,\tilde{K}$ epochs}
    \STATE{Calculate $L_{mse}^{[\tilde{k}]}=\frac{1}{n}\sqrt{\|\nabla P_{target} -\nabla {S_0}_{InitialNet}\|_2}$}
    \STATE{Update ${S_0}_{InitialNet}$ by gradient descent:$\Theta_{init}^{[\tilde{k}+1]} = \Theta_{init}^{[\tilde{k}]}-\eta_{init} \nabla_{\Theta_{init}} L_{mse}^{[\tilde{k}]}$}
    \ENDFOR
    \STATE{${S_0}^{[0]}_{NN}(\Theta) = {S_0}_{InitialNet}(\Theta_{init}^{[\tilde{K}]})$}
  \end{algorithmic}
\end{algorithm}

\section{Numerical Experiments}\label{sec:experiments}
We apply our method to Graph MFG with various graphs, including regular graphs (lattice and triangle grids) and inhomogeneous graphs, and different potentials.
In this work, we consider the edges on the graphs to be equally weighted, $G=(E,V,\mathbf{1})$, \textit{i.e.} $w_{ij}=w_{ji}=1, \forall (i,j)\in E.$


In our experiments, the fully-connected multi-layer-perceptron (MLP) is applied on regular graphs. For inhomogeneous graphs, a more sophisticated network designed for graphs with message passing and aggregation of the neighbors' information is preferable, thus we apply a graph neural network (GNN)~\cite{scarselli2008graph}, namely Graph Sample and Aggregate (GSAGE)~\cite{hamilton2017inductive}. 

\textit{\underbar{Multi-Layer Perceptron (MLP):}}
We choose the node coordinates, $X_i = (x_i,y_i)\in \mathbb{R}^2$ as the 2-dimensional input indices in this work and ${S_i}_0$ as the output:
${S_0}_{NN} (X_i) = S_0, i \in V.$ We use the 2-hidden-layer MLP: $${S_0}_{NN}(X_i) =  \mathbf{W}_{out}\cdot \sigma \left( \mathbf{W}_2 \cdot \sigma \left( \mathbf{W}_1 \cdot X_i + \mathbf{b}_1\right) + \mathbf{b}_2 \right) + \mathbf{b}_{out},$$ where $\sigma(\cdot) = \text{Tanh}(\cdot)$, and $\mathbf{W}$'s and $\mathbf{b}$'s are weights and bias. We choose the hidden neurons number to be 16, and the total parameters in the MLP is 337.
 This architecture is utilized in the experiments on regular graphs.

\textit{\underbar{Graph Sample and Aggregate (GSAGE):}}
For the GSAGE~\cite{hamilton2017inductive}, after the raw feature embedding of each node ($X_i \in \mathbb{R}^{2} \rightarrow \mathbf{h}_i^{(0)} \in \mathbb{R}^{16}$), there are two steps in each epoch, including neighbor nodes sampling and feature aggregation. 
As a GNN, the input includes indices $X_i$, $i$ in the node set $V$, and the edge set $E$, and the output is ${S_0}_i$, \textit{i.e.} ${S_0}_{GSAGE} (X_i, E) = {S_0}_i, i \in V.$
We use the GSAGE model with 3 layers that enable every node reaching to 3 hops of neighbor nodes, in which the feature flow is
$$\mathbf{h}_{i}^{(l)} = \sigma_{R} \left( \mathbf{W}_l \cdot \left[ \mathbf{h}^{(l-1)}_i \| \text{AGG}(\{\mathbf{h}^{(l-1)}_j, \forall j \in N(i)\}) \right] + \mathbf{b}_l \right)$$ at the $l$th layer, $l= 1,2, 3,$
followed by the output MLP layer ${S_0}_{GSAGE}(X_i,E) = \mathbf{W}_{out} \cdot \sigma(\mathbf{h}_{i}^{(3)})+\mathbf{b}_{out},$ where $\sigma_{R}(\cdot)=\text{RELU}(\cdot), \sigma(\cdot)=\text{Tanh}(\cdot)$, $\mathbf{h}_{i}^{(1)} = X_i, i \in V$. The $\text{AGG}$ here is the aggregation function (`MEAN' in our implementation) and $[ \cdot \| \cdot ]$ is concatenation operation. 
In this 3-layer GSAGE network, the number of parameters is 745.


\subsection{Test Examples on Lattice Graphs}\label{subsec:}

In these numerical tests, the graphs are lattice graphs embedded in a square $\mathcal{O}=[-3,3]\times[-3,3]$. The given initial density and the target terminal density follow 2-variate Gaussian distributions $\rho_G(\cdot,\mathbf{m},\mathbf{\Sigma}),$ in which mean $\mathbf{m} \in \mathbb{R}^2$ and covariance matrix $\mathbf{\Sigma} \in \mathbb{R}^{2 \times 2}.$ Specifically, given $\mu_0(x) = \rho_G (x,\mathbf{-1.2},0.2\cdot \mathbf{I}_2),$ the target density is $\mu_T(x) = \rho_G (x,\mathbf{1.2},0.2\cdot \mathbf{I}_2),$ where the bold number is the vector consisted of the same value, \textit{e.g.} $\mathbf{r}=r (\mathbf{e}_1+\mathbf{e}_2)$. The $\mathbf{e}_1$ and $\mathbf{e}_2$ are the first two standard basis vectors. The training settings are learning rate $lr=10^{-3}$ and training epochs 3000, and the settings remains unchanged in the subsequent experiments if not specified.

\subsubsection{Various Time Horizons and Step Sizes}\label{subsec:T-dt}

In this example, the graph $G$ is a $31\times31$ lattice graph, which has 961 nodes and 1860 undirected edges.
In the loss function in these experiments, the terminal energy is $\mathcal{G}_{L_1}$, and the coefficient ratio is $\lambda_K:\lambda_G=0.5:5000,$ where we denote the weight of the kinetic energy $\lambda_K=0.5$ based on the~\cref{def:iv-Potential MFG-d}.



\begin{table}[htbp]
	\footnotesize
	\begin{center}
		\begin{tabular}{|c|c|c|c|c|} \hline
			\bf $T$ & \bf $\Delta t$ & \bf $M$ & \bf Kinetic $\int_0^T \mathcal{K} \,dt$   & \bf Terminal $\mathcal{G}_{L_1}$\\ \hline
      25 & 0.125 & 200 & 23.1297 & 0.0229 \\ \hline
			25 & 0.25 & 100 & 23.3124 & 0.0184 \\ \hline
			25 & 0.5 & 50 & 23.6192 & 0.0332 \\ \hline
      50 & 0.25 & 200 & 11.5819 & 0.0134 \\ \hline
      50 & 0.5 & 100 & 11.6432 & 0.0126 \\ \hline
      50 & 1.0 & 50 & 11.8686 & 0.0217  \\ \hline
			100 & 0.25 & 400 & 5.7853 & 0.0113  \\ \hline
			100 & 0.5 & 200 & 5.7933 & 0.0125 \\ \hline
			100 & 1.0 & 100 & 5.8181  & 0.0191 \\ \hline
		\end{tabular}
		\caption{Comparisons of the total kinetic energy $\int_0^T \mathcal{K} \,dt$ and the terminal energy $\mathcal{G}_{L_1}$ on the lattice graph $G$ with different time horizons $T$, time interval $\Delta t$ and steps number $M$. These results are the average of the last 100 epochs of the total 3000 epochs.}\label{tab:compare-ss}
	\end{center}
\end{table}

In~\cref{{tab:compare-ss}}, we examine the appropriate time horizons $T$, time interval $\Delta t$. In these experiments, all the crowd successfully arrive the target terminal density $\mu_T$ after the training convergence. When $T$ is doubled, the total kinetic energy is about halved as shown in~\cref{tab:compare-ss}. It can be understood as follows. When the distance is fixed, the velocity should be inversely propotional to the total time horizon, \textit{i.e.} $v \propto \frac{1}{T}$, where $M=\frac{T}{\Delta t}$. Therefore, when $T$ is doubled, velocity should be halved, and accordingly the total kinetic energy $\int_0^T \mathcal{K} dt = \sum_{m=0}^M \sum_{(i,j)\in E} \frac{1}{2}\|v_{ij}\|^2 \theta_{ij}(\rho) \Delta t$ should be about halved. Besides, we can also see in~\cref{tab:compare-ss} that larger $T$ enables smaller terminal errors. 


\subsubsection{Different Graph Sizes}
We consider different grid sizes of the lattice graphs, from $11\times 11$ to $41\times 41$. 
In the loss function, the terminal energy is $\mathcal{G}_{L_1}$, and the coefficient ratio is $\lambda_K:\lambda_G=0.5:5000.$
These experiments have the uniform settings that time step $\Delta t=0.5$ and steps number $M=200$. 



\begin{table}[htbp]
	\footnotesize
	\begin{center}
		\begin{tabular}{|c|c|c|} \hline
			\bf Grid Size & \bf Kinetic $\int_0^T \bar{\mathcal{K}} \,dt$ & \bf Terminal $\bar{\mathcal{G}}_{L_1}$  \\ \hline
			$11\times11$ & $7.5701 \times 10^{-3}$ & $4.5386 \times 10^{-4}$  \\\hline
			$21\times21$ & $6.1428\times 10^{-3}$  & $6.3765 \times 10^{-5}$ \\ \hline
			$31\times31$ & \textbf{$6.0284 \times 10^{-3}$} & \textbf{$1.2957 \times 10^{-5}$} \\ \hline
			$41\times41$ & $6.1112 \times 10^{-3}$ & \textbf{$7.6677\times 10^{-6}$} \\ \hline
		\end{tabular}
		\caption{Results of lattice graphs with different grid sizes. The table displays node-average values of the total kinetic energy $\int_0^T \bar{\mathcal{K}}\,dt$ and the terminal energy $\bar{\mathcal{G}}_{L_1}$.}\label{tab:compare-gs}
	\end{center}
\end{table}
The node-average of total kinetic energy $\int_0^T \bar{\mathcal{K}} \,dt$ and the terminal energy  $\bar{\mathcal{G}}_{L_1}=\frac{1}{n}\|\rho(T)-\mu_T\|_1$ are compared in the~\cref{tab:compare-gs}. 
\Cref{tab:compare-gs} shows that the node-average terminal energy decreases as the number of nodes increase. The results in the table are the costs evaluated on the best model after the training convergence. 

\subsubsection{Different Kinetic Ratios and Different Terminal Energies}
The graph settings are the same as those in the example in~\cref{subsec:T-dt}.
These experiments use the uniform settings that $\Delta t=0.5$ and steps number $100$. In the loss function, different ratios of kinetic to terminal energies  of $L_1$ distance $\mathcal{G}_{L_1}$ and KL divergence $\mathcal{G}_{KL}$ are compared, \textit{i.e.} $\lambda_K:\lambda_G.$

\begin{table}[htbp]
	\footnotesize
	\begin{center}
    \begin{minipage}{0.45\textwidth}
    \vspace{0pt}
    \centering
		\resizebox{\textwidth}{!}{\begin{tabular}{|c|c|c|}
      \hline
			$\lambda_k:\lambda_g$  & \bf Kinetic $\int_0^T {\mathcal{K}} \,dt$ & \bf Terminal $\mathcal{G}_{L_1}$ \\ \hline
      $0:50$ & 11.6431 & 0.0200 \\ \hline
      $0:5000$ & 11.6321 & 0.0144 \\ \hline
      $0.5:5000$ &  11.6432  & 0.0126 \\ \hline
      $0.5:500$ & 11.6376 & 0.0161 \\ \hline
      $0.5:50$ & 11.5931 & 0.0207 \\ \hline
		\end{tabular}}\subcaption{Terminal energy $\mathcal{G}_{L_1}$.}\label{tab:compare-losses}
    \end{minipage}%
    \hfill
    \begin{minipage}{0.45\textwidth}
    \vspace{0pt}
    \centering
    \resizebox{\textwidth}{!}{\begin{tabular}{|c|c|c|c|}
    \hline
    $\lambda_k:\lambda_g$  & \bf Kinetic $\int_0^T {\mathcal{K}} \,dt$ & \bf Terminal $\mathcal{G}_{KL}$ \\ \hline
    $0:50$ & 11.7042 & $4.0322\times 10^{-4}$ \\ \hline
    $0:5000$ & 11.6840 & $2.5029\times 10^{-4}$  \\ \hline
    $0.5:5000$ & 11.6741 & $3.5065\times 10^{-4}$  \\ \hline
    $0.5:500$ & 11.6372 & $4.5081\times 10^{-4}$ \\ \hline
    $0.5:50$ & 11.5701 & $9.1995\times 10^{-4}$ \\ \hline
    \end{tabular}}\subcaption{Terminal energy $\mathcal{G}_{\mathcal{KL}}$.}\label{tab:compare-losses-kld}
    \end{minipage}
		\caption{Comparisons between the experiments on different ratios of kinetic energy and terminal energy ($\lambda_K:\lambda_G$), and different terminal energies $\mathcal{G}_{L_1}$ and $\mathcal{G}_{\mathcal{KL}}$.}\label{tab:compare-kin-ratio}
	\end{center}
\end{table}



%

\Cref{tab:compare-losses,tab:compare-losses-kld} display the result comparisons of different ratios of the energies, and different terminal energies $\mathcal{G}_{L_1}$ and $\mathcal{G}_{KL}$. Noticed that when the kinetic energy has coefficient 0, the kinetic and terminal energies also converge to a small value. This is because the critical point conditions~\cref{eq:fw-hj-ce} ensures the lowest kinetic energy. When the objective function is convex, because the solution is unique, both of the paths with and without kinetic energy are the solution to the Graph MFG. 
It can be seen that a good result is obtained under the condition $\lambda_K:\lambda_G=0.5:5000$, and both $\mathcal{G}_{L_1}$ and $\mathcal{G}_{KL}$ give a good constraint at the target density. Therefore, the ratio $\lambda_K:\lambda_G=0.5:5000$ and terminal condition $\mathcal{G}_{L_1}$ are chosen defaultly for all Graph MFG experiments in~\cref{experiments:GMFG}.

\subsection{Graph MFG Experiments on Different Graphs}\label{experiments:GMFG}

\subsubsection{Graph MFG on Regular Graphs}

\begin{example}{Lattice Grids with Square Boundary}

In this example, the graph $G$ is embedded in the same spatial domain as that in~\cref{subsec:T-dt}, except it has 1681 nodes and 3280 undirected edges.
Initial and target densities of the agents are also the same as the previous examples.
The experimental settings are $T=50$, $\Delta t=0.5$ and $M=100$. 

\Cref{fig:re-lattice} shows the dynamics of the population and value function. The obtained path of the crowd moves straightly from the given initial density to the target, without significant change of size. In this result, the total kinetic energy is 10.2583, terminal energy is 0.0139 and $\frac{1}{n}\|\rho(T)-\mu_T\|_1=8.3109*10^{-6}$. 



\begin{figure*}[htbp] 
  \centering
  \begin{minipage}[b]{0.99\textwidth} 
  {
    \begin{minipage}[b]{0.14\linewidth} 
      \centering
      \includegraphics[width=\linewidth]{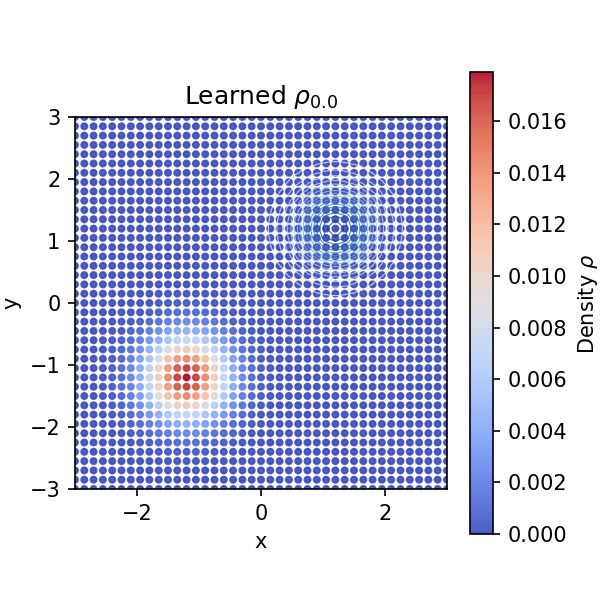}
      \includegraphics[width=\linewidth]{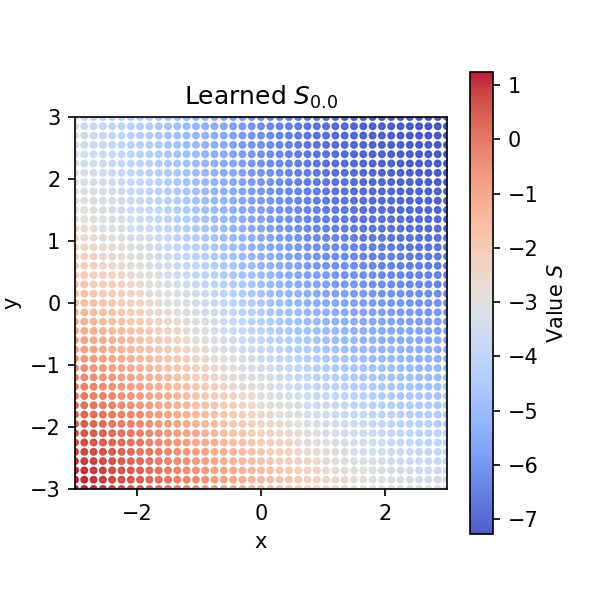}\subcaption{$t=0$}
    \end{minipage}
  }
   {
    \begin{minipage}[b]{0.14\linewidth}
      \centering
      \includegraphics[width=\linewidth]{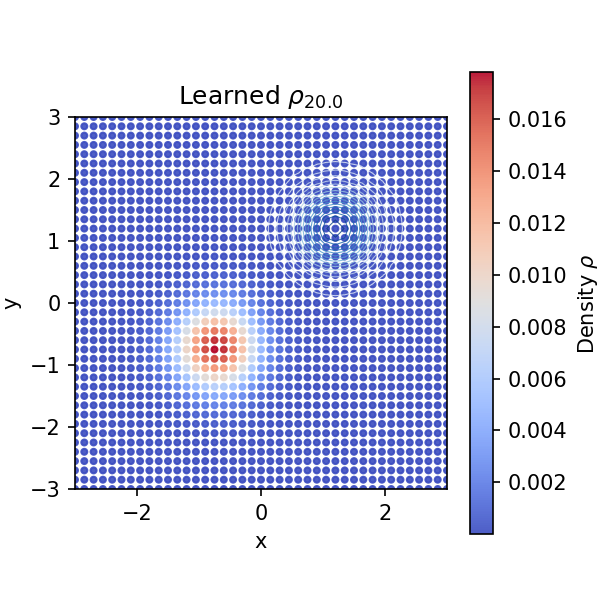}
      \includegraphics[width=\linewidth]{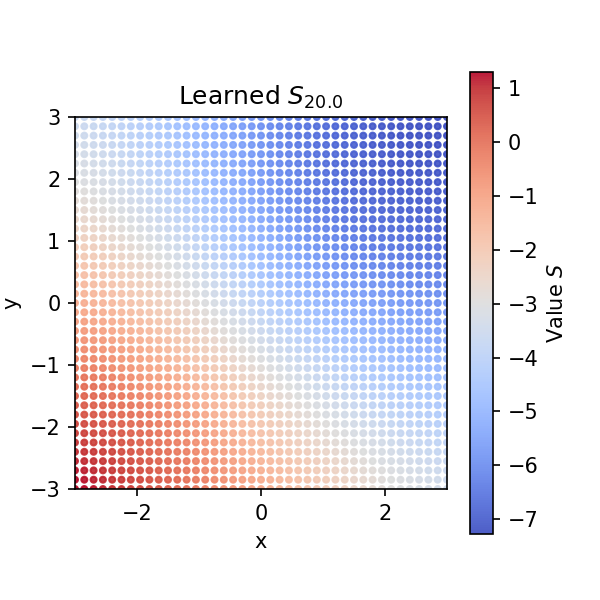}\subcaption{$t=20$}
    \end{minipage}
  }
   {
    \begin{minipage}[b]{0.14\linewidth}
      \centering
      \includegraphics[width=\linewidth]{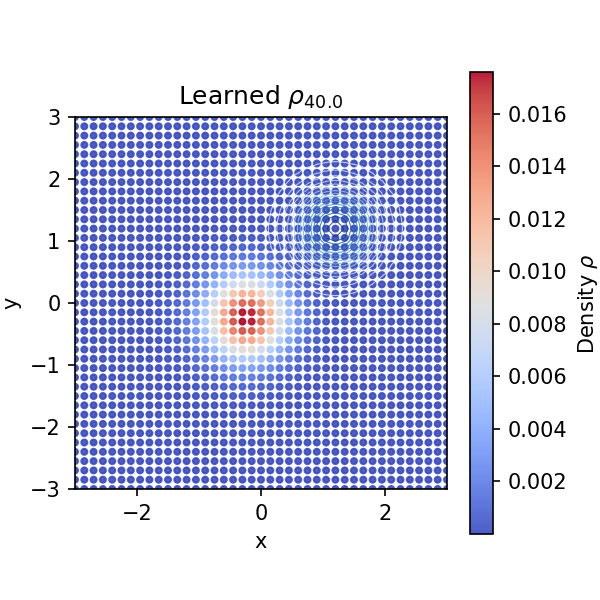}
      \includegraphics[width=\linewidth]{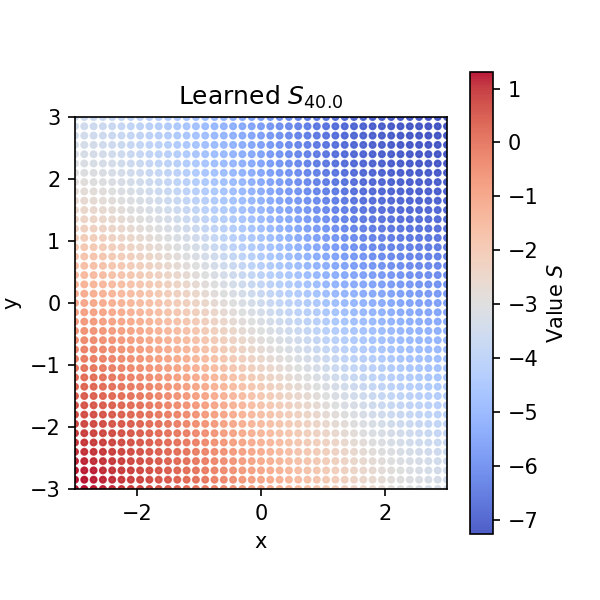}\subcaption{$t=40$}
    \end{minipage}
  }
    {
    \begin{minipage}[b]{0.14\linewidth}
      \centering
      \includegraphics[width=\linewidth]{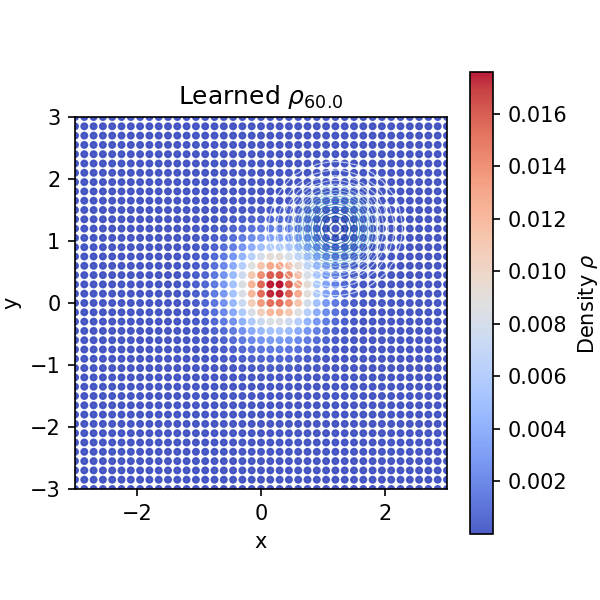}
      \includegraphics[width=\linewidth]{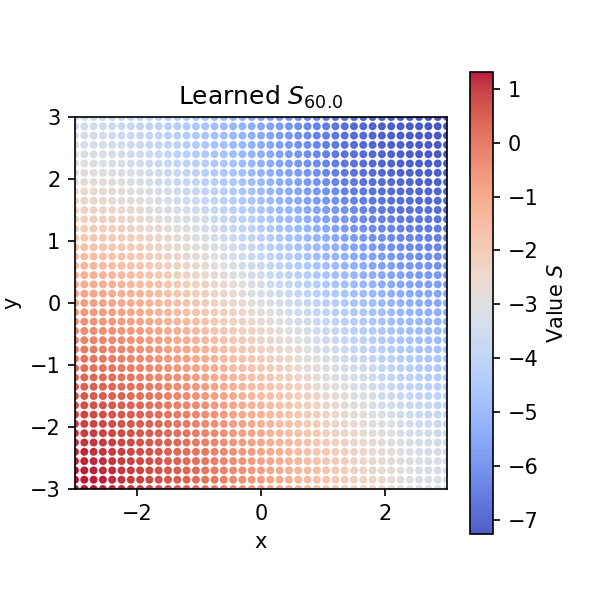}\subcaption{$t=60$}
    \end{minipage}
  }
    {
    \begin{minipage}[b]{0.14\linewidth}
      \centering
      \includegraphics[width=\linewidth]{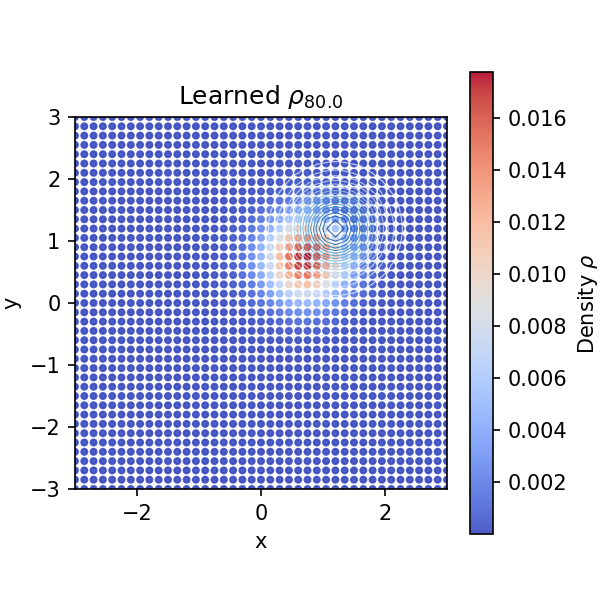}
      \includegraphics[width=\linewidth]{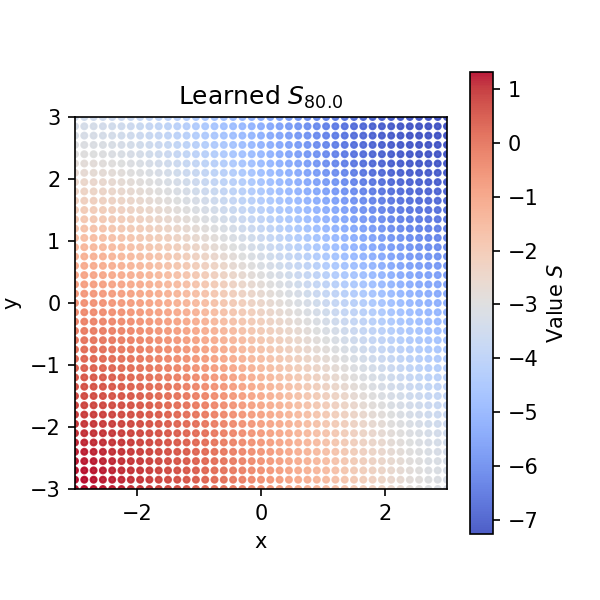}\subcaption{$t=80$}
    \end{minipage}
  }
     {
    \begin{minipage}[b]{0.14\linewidth}
      \centering
      \includegraphics[width=\linewidth]{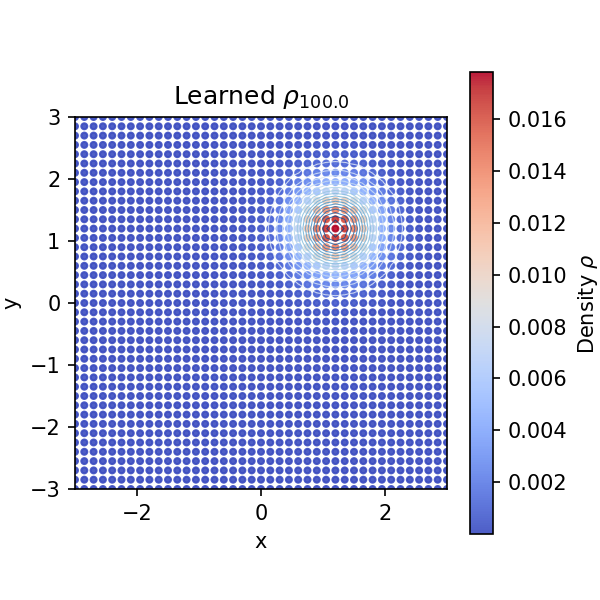}
      \includegraphics[width=\linewidth]{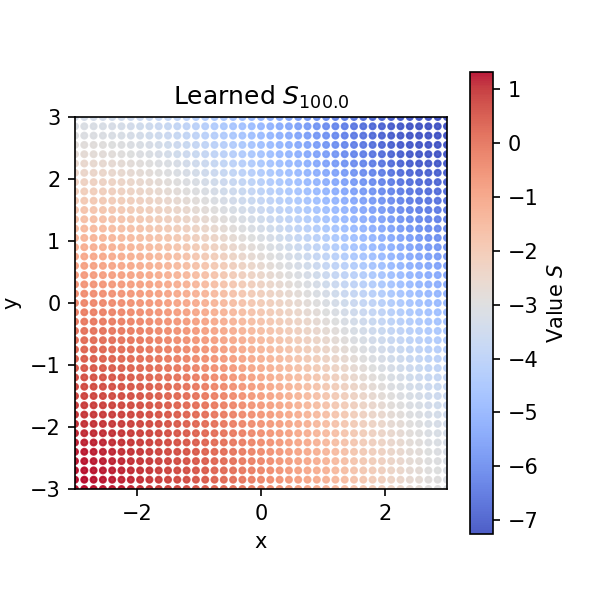}\subcaption{$t=100$}
    \end{minipage}
    }
  \end{minipage}
  \caption{Results of the dynamics of population density $\rho$ (upper panel) and value function $S$ (lower panel) on the regular graph with lattice grids.}\label{fig:re-lattice}
\end{figure*}

\end{example}

\begin{example}{Triangle Grids with Irregular Boundary}


  In this example, the graph $G$ is a triangle graph embedded in $\mathcal{\hat{O}}=x^2+y^2<3$, with a rectangle vacancy $\mathcal{R}=[-3,-1.5] \times [-1,1]$ and a circle hole inside $\tilde{\mathcal{O}}=x^2+y^2<0.8$, thus the spatial domain is $\mathcal{O}=\hat{O}\backslash \{\mathcal{R}\bigcup \mathcal{\tilde{O}} \}$. Here, $G$ has 463 nodes and 1272 undirected edges, which consists of the uniform triangle meshes and irregular boundary shape.
  The initial density obeys the Gaussian mixture distribution, $\mu_0 (x) =\frac{1}{3} \sum_{i=1}^3 \rho_G(x,\mathbf{m}_i,0.15\cdot \mathbf{I}_2),$ where $\mathbf{m}_i = \mathbf{p}_0 + r(\sin(\theta_n)\mathbf{e}_1+\cos(\theta_i)\mathbf{e}_2), \mathbf{p}_0=\mathbf{0}, r=2, \theta_n=\frac{2\pi}{3} i, i=1,2,3.$ The target terminal density obeys the Laplacian distribution, $\mu_T(x) = \text{exp}(-a_0|x_0-b_0|-a_1|x_1-b_1|),$ in which $a_0=a_1=2, b_0=b_1=2.$
  The experimental settings are $T=20$, $\Delta t=0.5$ and $M=40$. 
  The training epochs is 30000. 

  \cref{fig:re-triangle} shows the dynamics of the population density and value function. There are three groups of crowds at the beginning. The upper two groups go straightly to the target, while the left-bottom one moves right and merge with the other groups, because there are vacancies on ths graph that the crowd can not move through them. In this result, the total kinetic energy is 7.1114, terminal energy is 0.4270, and $\frac{1}{n}\|\rho(T)-\mu_T\|_1=0.0009$. 
  
  \begin{figure*}[htbp] 
    \centering
    \begin{minipage}[b]{0.99\textwidth} 
    {
      \begin{minipage}[b]{0.14\linewidth} 
        \centering
        \includegraphics[width=\linewidth]{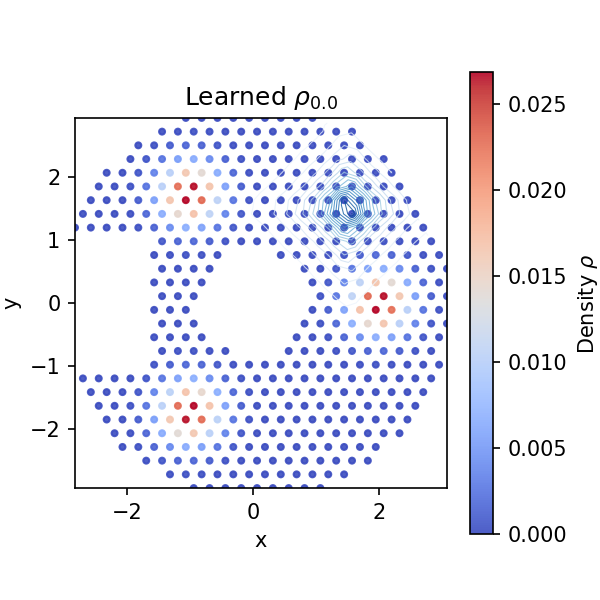}
        \includegraphics[width=\linewidth]{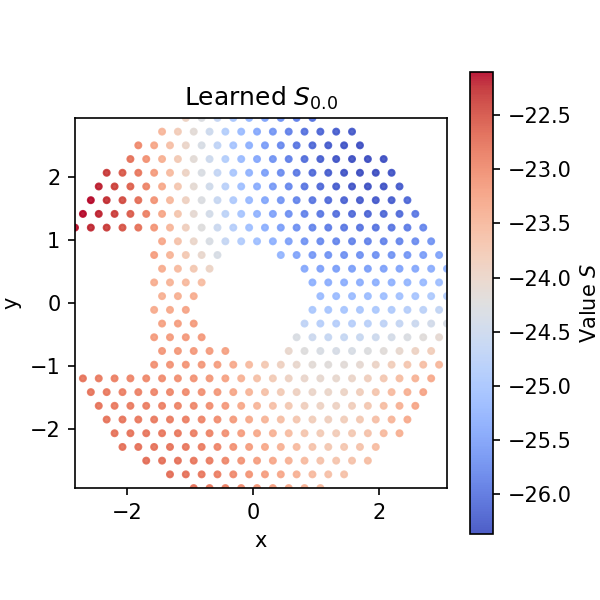}\subcaption{$t=0$}
      \end{minipage}
    }
     {
      \begin{minipage}[b]{0.14\linewidth}
        \centering
        \includegraphics[width=\linewidth]{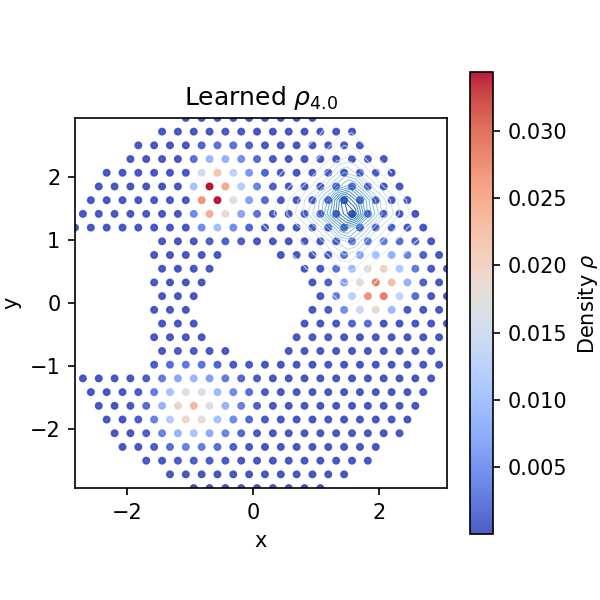}
        \includegraphics[width=\linewidth]{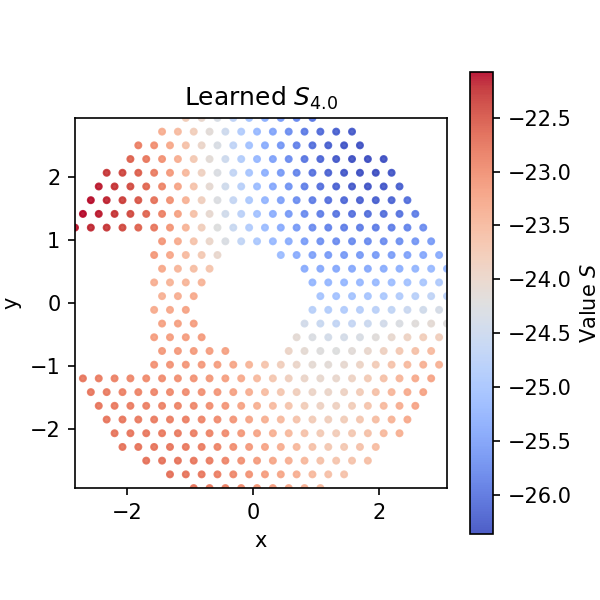}\subcaption{$t=4$}
      \end{minipage}
    }
     {
      \begin{minipage}[b]{0.14\linewidth}
        \centering
        \includegraphics[width=\linewidth]{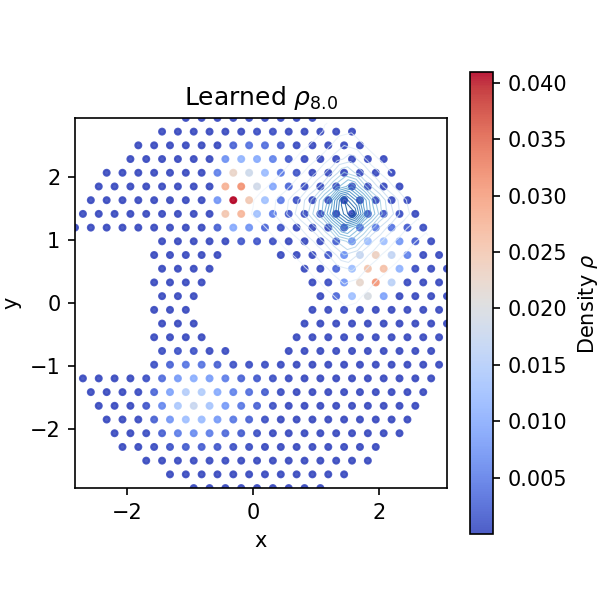}
        \includegraphics[width=\linewidth]{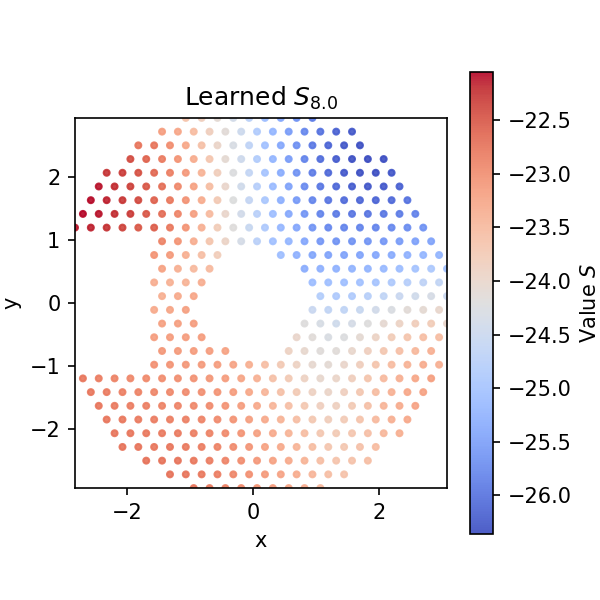}\subcaption{$t=8$}
      \end{minipage}
    }
      {
      \begin{minipage}[b]{0.14\linewidth}
        \centering
        \includegraphics[width=\linewidth]{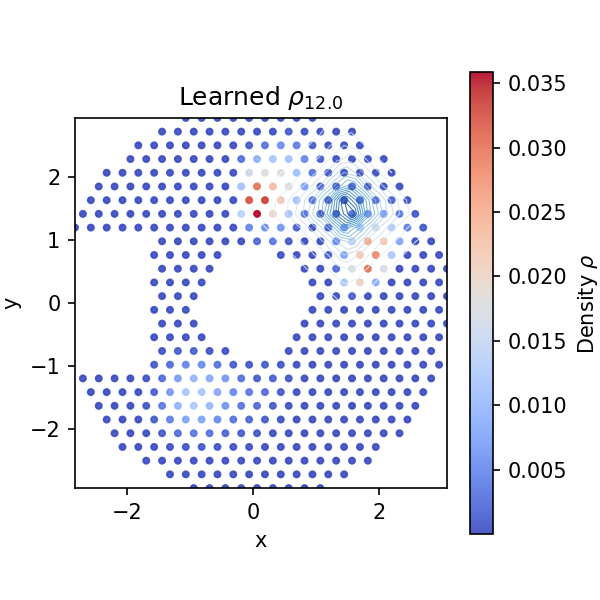}
        \includegraphics[width=\linewidth]{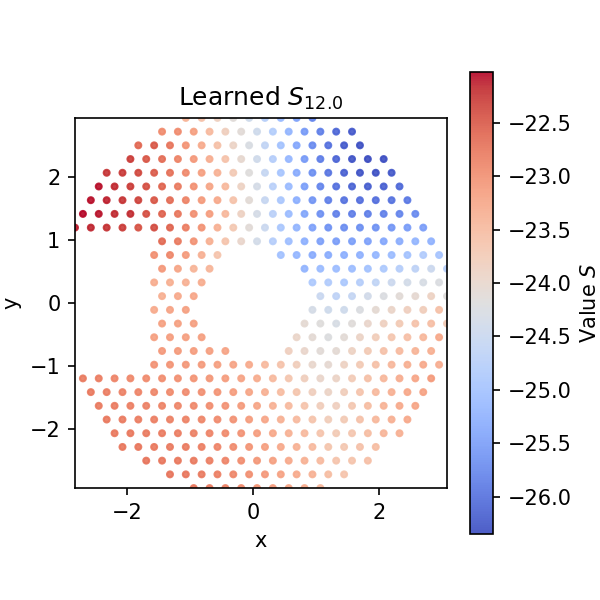}\subcaption{$t=12$}
      \end{minipage}
    }
      {
      \begin{minipage}[b]{0.14\linewidth}
        \centering
        \includegraphics[width=\linewidth]{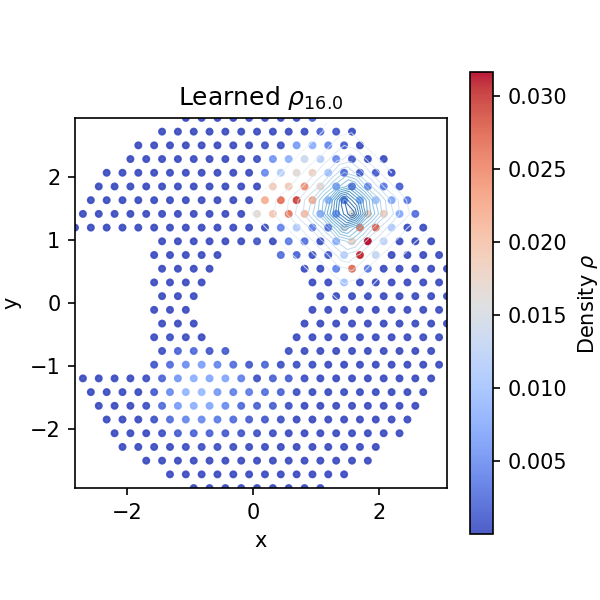}
        \includegraphics[width=\linewidth]{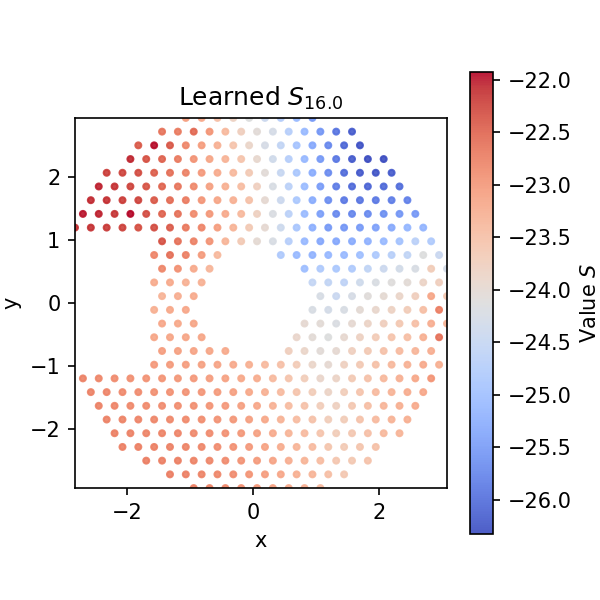}\subcaption{$t=16$}
      \end{minipage}
      }
      {
      \begin{minipage}[b]{0.14\linewidth}
        \centering
        \includegraphics[width=\linewidth]{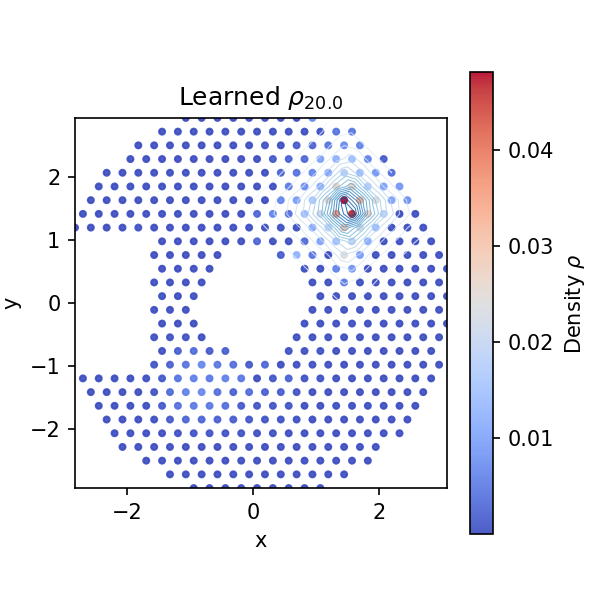}
        \includegraphics[width=\linewidth]{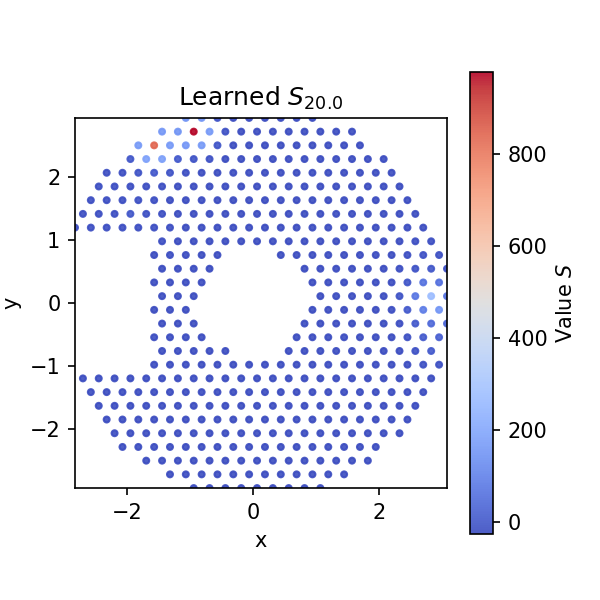}\subcaption{$t=20$}
      \end{minipage}
    }
    
    \end{minipage}
    \caption{Result of the dynamics of population density $\rho$ (upper panel) and value function $S$ (lower panel) on the regular graph with triangle grids and specific-shaped boundary.}\label{fig:re-triangle}
  \end{figure*}
  
  \end{example}
  
\begin{example}{Triangle Grids and Linear Potential}

In this example, the graph $G$ is a triangle graph embedded in a circle $\mathcal{O}=x^2+y^2<3$. It has 563 nodes and 1604 undirected edges.
The initial density of the crowd is $\mu_0(x) = \rho_G (x,\mathbf{-1},0.2\cdot \mathbf{I}_2),$ and the target is $\mu_T(x) = \rho_G (x,\mathbf{1},0.2\cdot \mathbf{I}_2).$ 
The linear potential $\mathcal{V}$ in this experiment is a smoothed uniform-like distribution that has a radius near $0.5$, playing the role of landscape preference potential as~\cite{lin2021alternating}. The landscape potential here is $\mathcal{V}(x) = 50\times \frac{1}{2}\left[ 1 + \text{Tanh}\left( \eta \times \left( R - \|x - \mathbf{c}\|_2 \right) \right) \right],$
where $R=0.5$ and $\mathbf{c}=0.5\mathbf{e}_1+0.5\mathbf{e}_2$ are the obstacle radius and center, and $\eta=5$ is the smoothness factor. 
We apply 3000 epochs in the warm up training, with the linear combination parameter $\alpha=0.6$ in~\cref{eq:warmup}. In the loss function of the main training, the coefficient ratio is $\lambda_K:\lambda_V:\lambda_G=0.5:0.01:5000.$
The experimental settings are $T=20$, $\Delta t=0.5$ and $M=40$, and the training epochs is 5000. We observe that it stabilizes the training and the loss decreases faster and converges to a better result. 

\Cref{fig:re-linear} shows the dynamics of the population density and value function. The crowd splits at the center and aggregates at the end to match the target, because the potential at the center is high. In the result, the total kinetic energy is 13.9315, linear potential of the dynamics is 10.6170, terminal energy is 0.2329, and $\frac{1}{n}\|\rho(T)-\mu_T\|_1=0.0004$. 
  
  \begin{figure*}[htbp] 
    \centering
    \begin{minipage}[b]{0.99\textwidth} 
    {
      \begin{minipage}[b]{0.15\linewidth} 
        \centering
        \includegraphics[width=\linewidth]{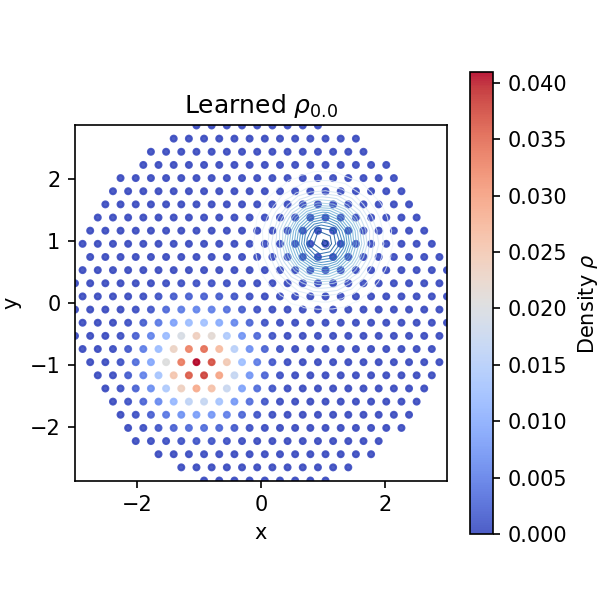}
        \includegraphics[width=\linewidth]{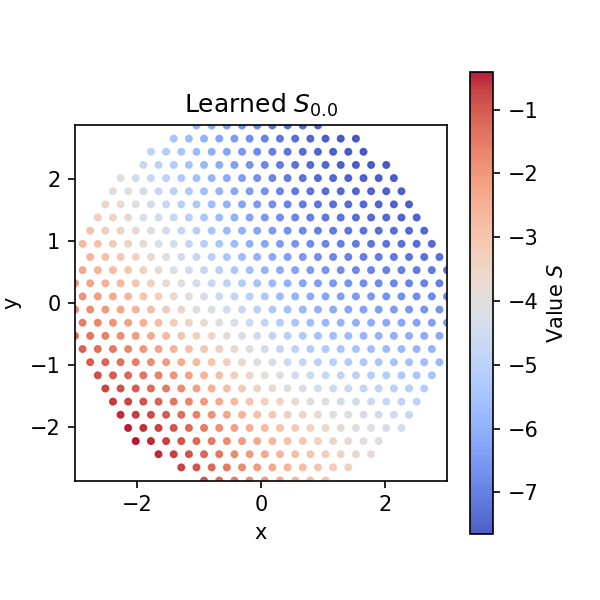}
        \subcaption{$t=0$}
      \end{minipage}
    }
     {
      \begin{minipage}[b]{0.15\linewidth}
        \centering
        \includegraphics[width=\linewidth]{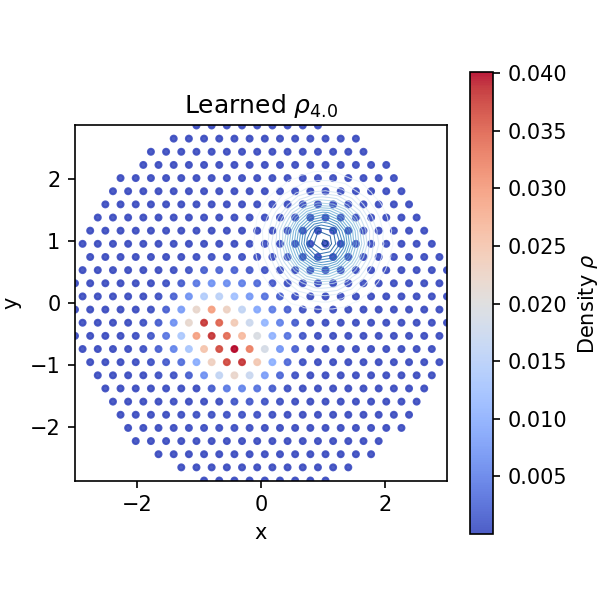}
        \includegraphics[width=\linewidth]{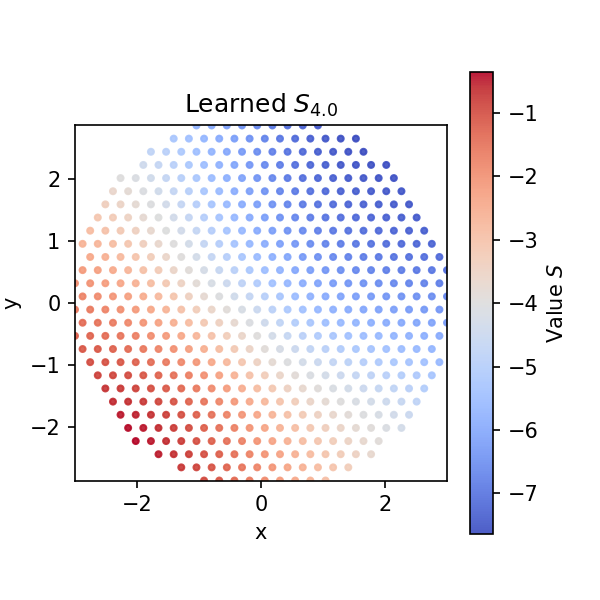}
        \subcaption{$t=4$}
      \end{minipage}
    }
     {
      \begin{minipage}[b]{0.15\linewidth}
        \centering
        \includegraphics[width=\linewidth]{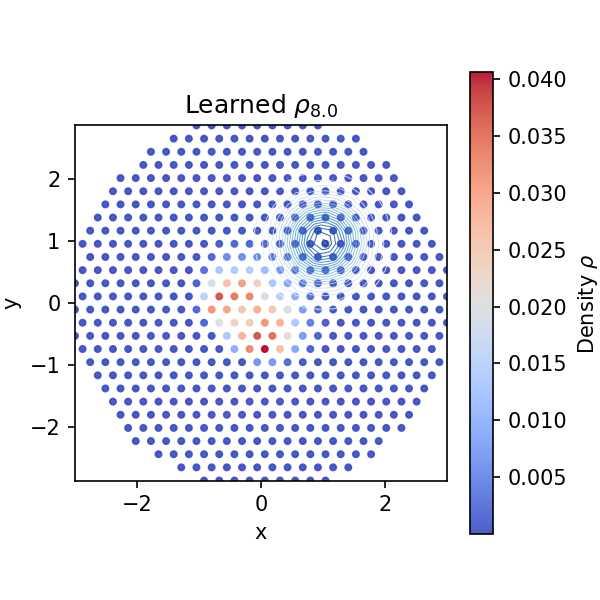}
        \includegraphics[width=\linewidth]{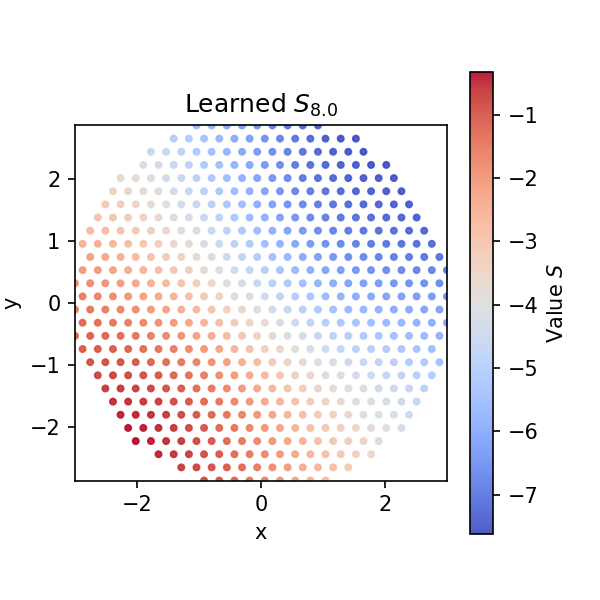}\subcaption{$t=8$}
      \end{minipage}
    }
      {
      \begin{minipage}[b]{0.15\linewidth}
        \centering
        \includegraphics[width=\linewidth]{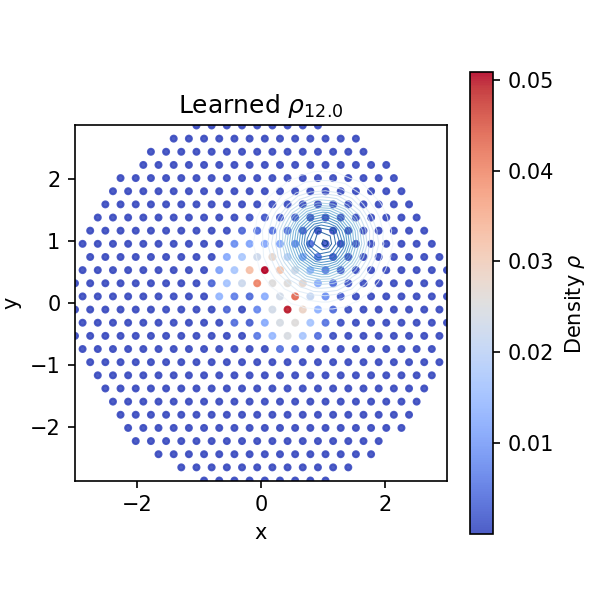}
        \includegraphics[width=\linewidth]{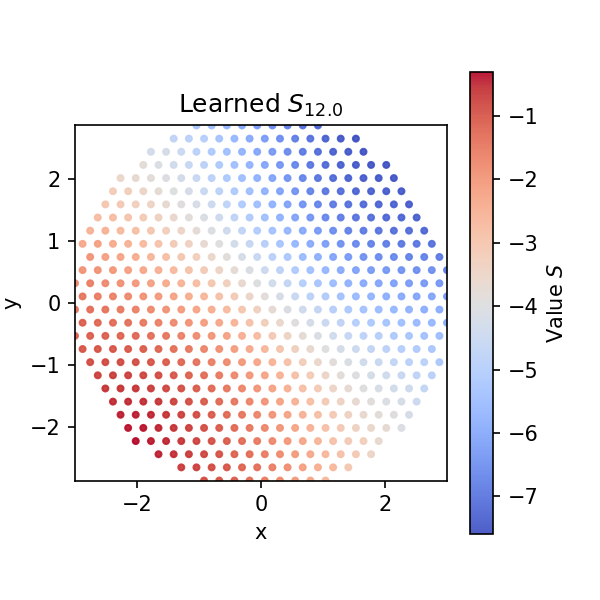}\subcaption{$t=12$}
      \end{minipage}
    }
      {
      \begin{minipage}[b]{0.15\linewidth}
        \centering
        \includegraphics[width=\linewidth]{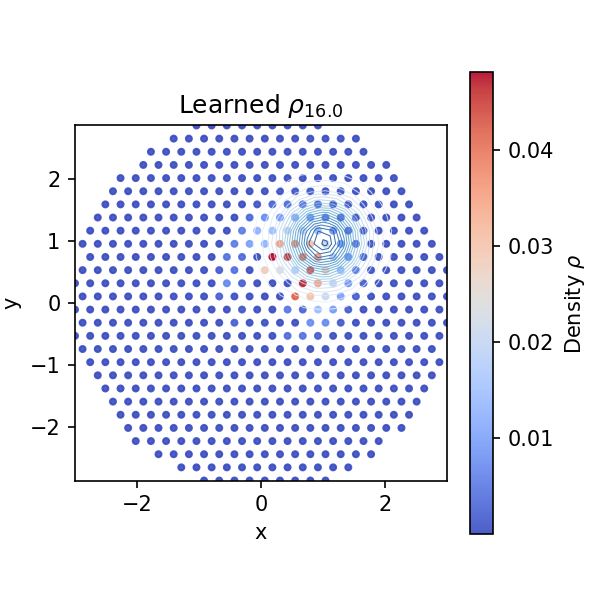}
        \includegraphics[width=\linewidth]{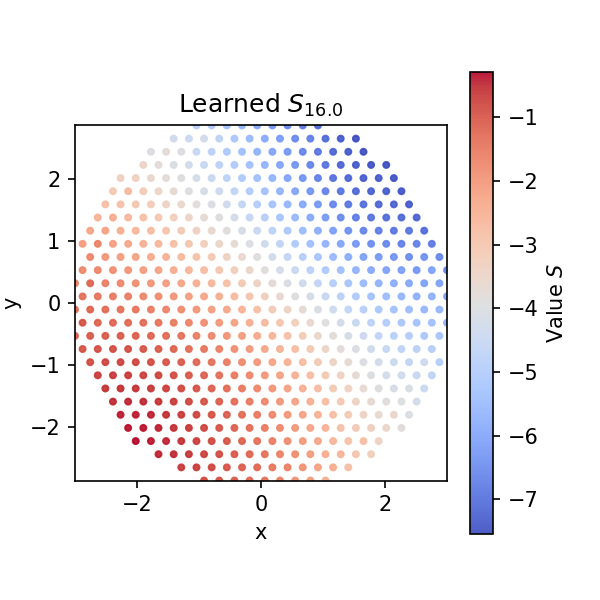}\subcaption{$t=16$}
      \end{minipage}
      \begin{minipage}[b]{0.15\linewidth}
       \centering
       \includegraphics[width=\linewidth]{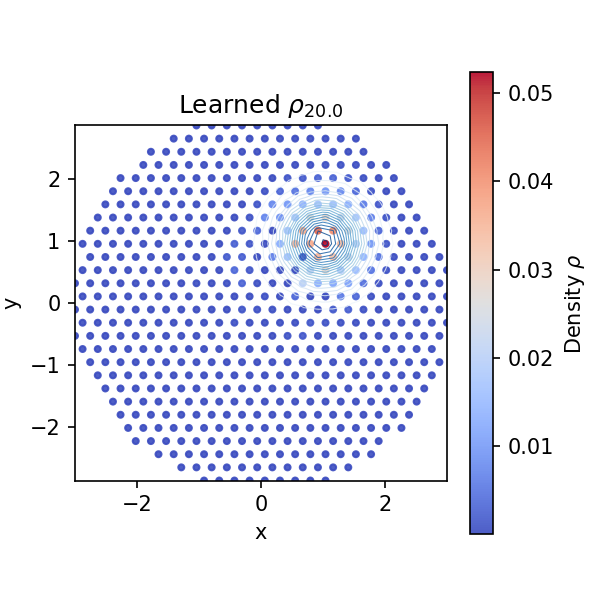}
       \includegraphics[width=\linewidth]{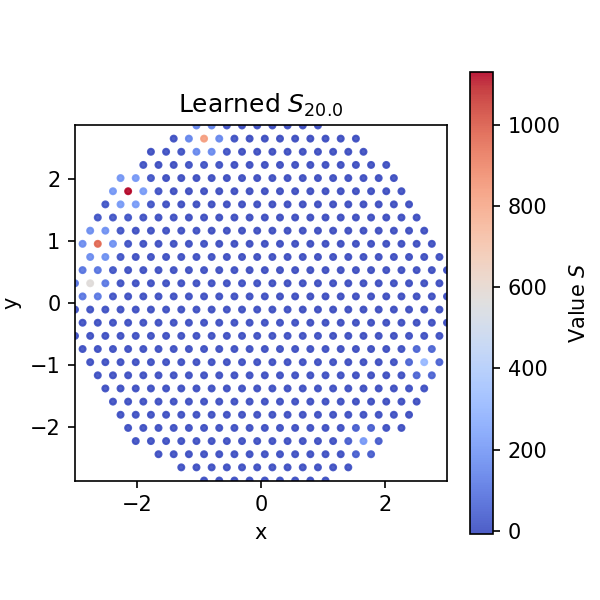}\subcaption{$t=20$}
      \end{minipage}
    }
    \end{minipage}
    \caption{Result of the dynamics of population density $\rho$ (upper panel) and value function $S$ (lower panel) with linear potential on the regular graph with triangle grids and circle boundary.}\label{fig:re-linear}
  \end{figure*}
  
  \end{example}

\begin{example}{Triangle Grids and Gaussian Interaction Potential}

In this example, the graph $G$ is a triangle graph embedded in a circle $\mathcal{O}=x^2+y^2<3$. It has 396 nodes and 1115 undirected edges.
The initial density of the crowd is $\mu_0(x) = \rho_G (x,\mathbf{-0.8},0.2\cdot \mathbf{I}_2),$ and the target is $ \mu_T(x) = \rho_G (x,\mathbf{0.8},0.2\cdot \mathbf{I}_2).$ 
In the loss function, the discrete Gaussian interaction potential is $\mathcal{W}(\rho):=\frac{1}{2} \sum_{j=1}^n \sum_{i=1}^n e^{-\frac{1}{2}\|X_i-Y_j\|^2} \rho_i \rho_j,$ and the coefficient ratio is $\lambda_K:\lambda_W:\lambda_G=0.5:0.1:5000.$
The experimental settings are $T=20$, $\Delta t=0.2$ and $M=100$.

\Cref{fig:re-Gaussian} shows the dynamics of the population density and value function. In this case, there is a Gaussian interaction potential that penalize aggregation. The croad spreads out from the center and successfully arrive to the target. 
In the result, the total kinetic energy is 4.1204, the interaction cost is 2.4693, terminal energy is 0.0239, and $\frac{1}{n}\|\rho(T)-\mu_T\|_1=6.0285\times 10^{-5}$. 


\begin{figure*}[htbp] 
  \centering
  \begin{minipage}[b]{0.99\textwidth} 
  {
    \begin{minipage}[b]{0.15\linewidth} 
      \centering
      \includegraphics[width=\linewidth]{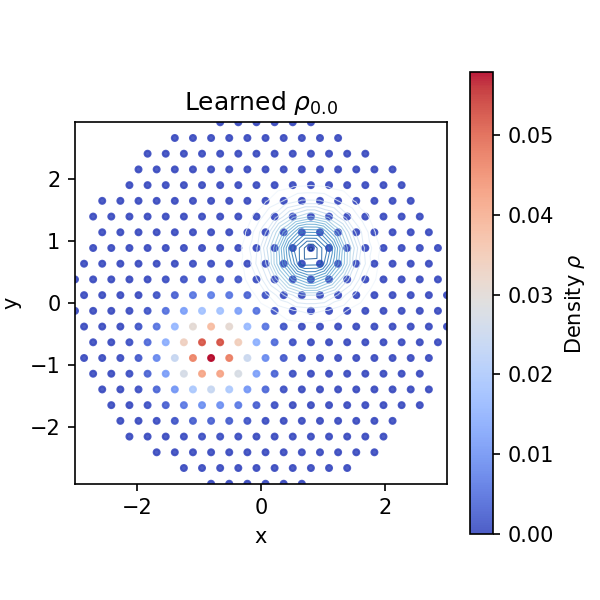}
      \includegraphics[width=\linewidth]{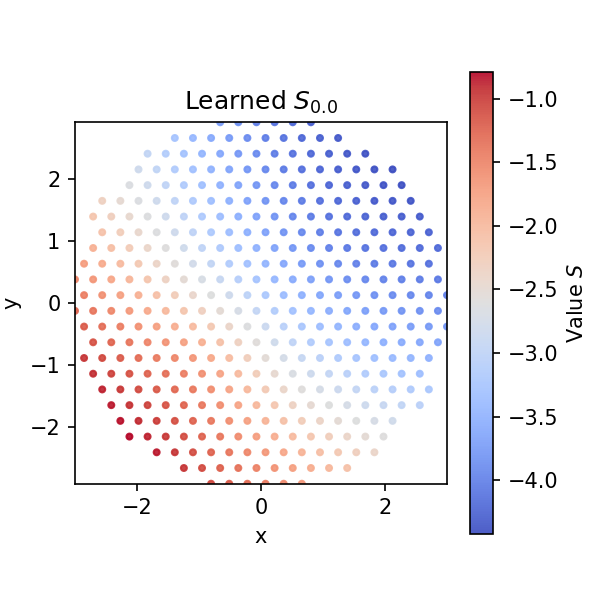}\subcaption{$t=0$}
    \end{minipage}
  }
   {
    \begin{minipage}[b]{0.15\linewidth}
      \centering
      \includegraphics[width=\linewidth]{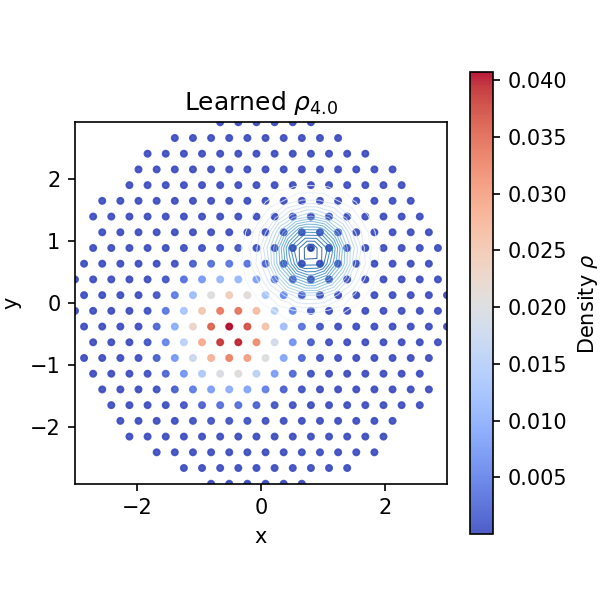}
      \includegraphics[width=\linewidth]{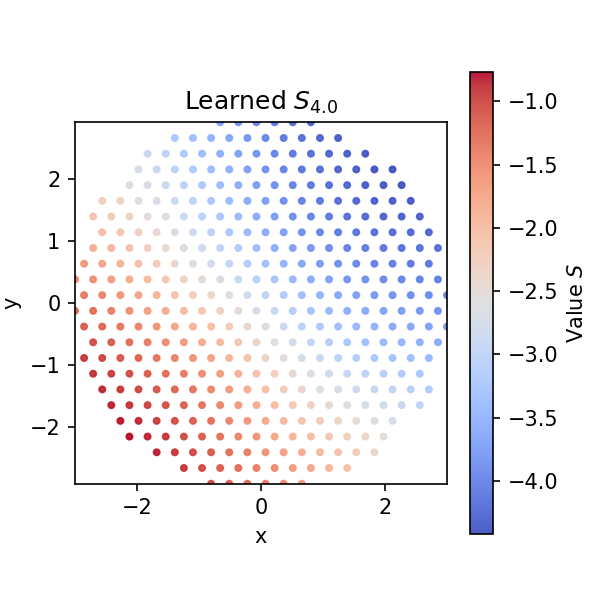}\subcaption{$t=4$}
    \end{minipage}
  }
   {
    \begin{minipage}[b]{0.15\linewidth}
      \centering
      \includegraphics[width=\linewidth]{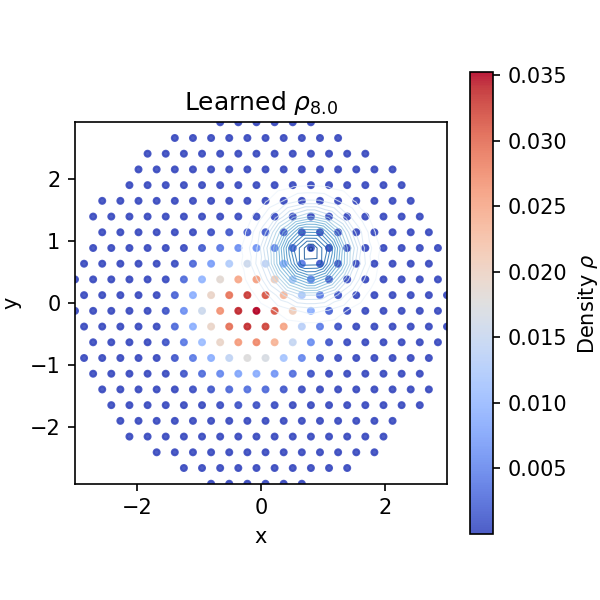}
      \includegraphics[width=\linewidth]{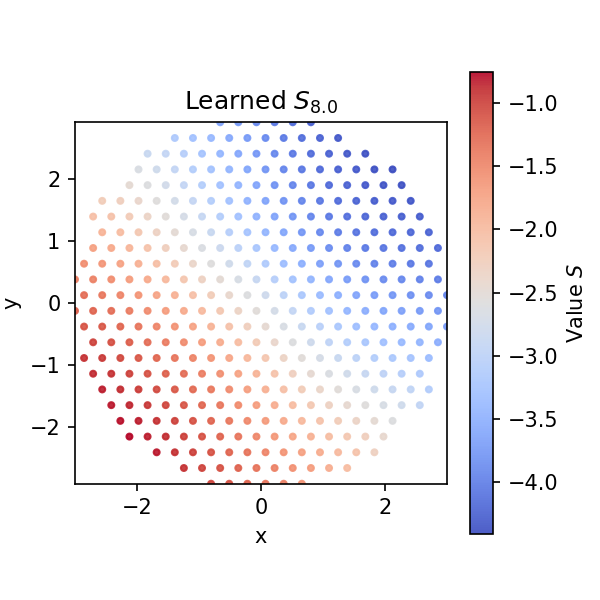}\subcaption{$t=8$}
    \end{minipage}
  }
    {
    \begin{minipage}[b]{0.15\linewidth}
      \centering
      \includegraphics[width=\linewidth]{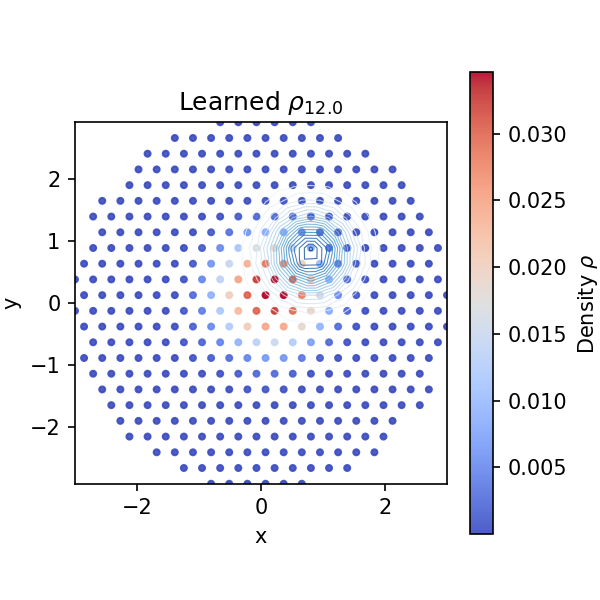}
      \includegraphics[width=\linewidth]{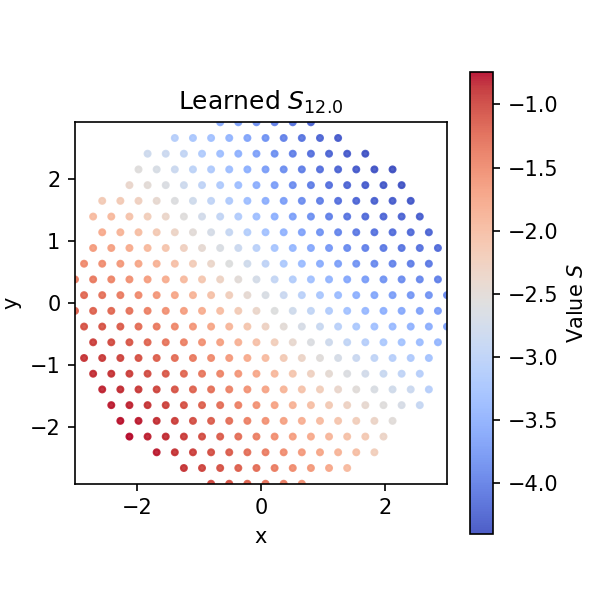}\subcaption{$t=12$}
    \end{minipage}
  }
    {
    \begin{minipage}[b]{0.15\linewidth}
      \centering
      \includegraphics[width=\linewidth]{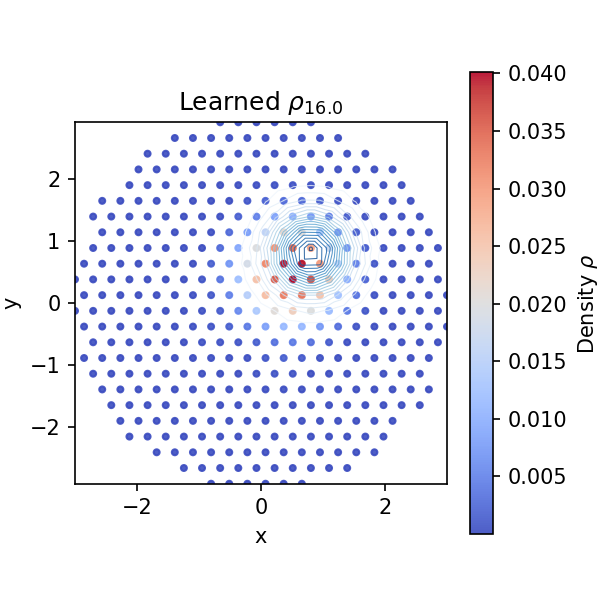}
      \includegraphics[width=\linewidth]{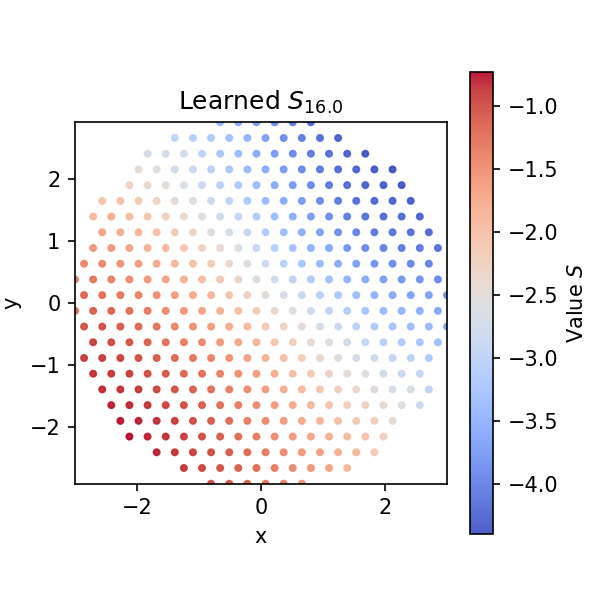}\subcaption{$t=16$}
    \end{minipage}
    \begin{minipage}[b]{0.15\linewidth}
      \centering
      \includegraphics[width=\linewidth]{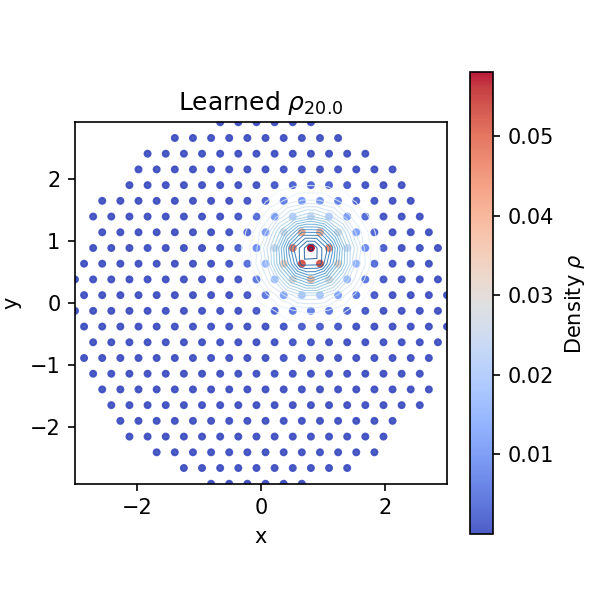}
      \includegraphics[width=\linewidth]{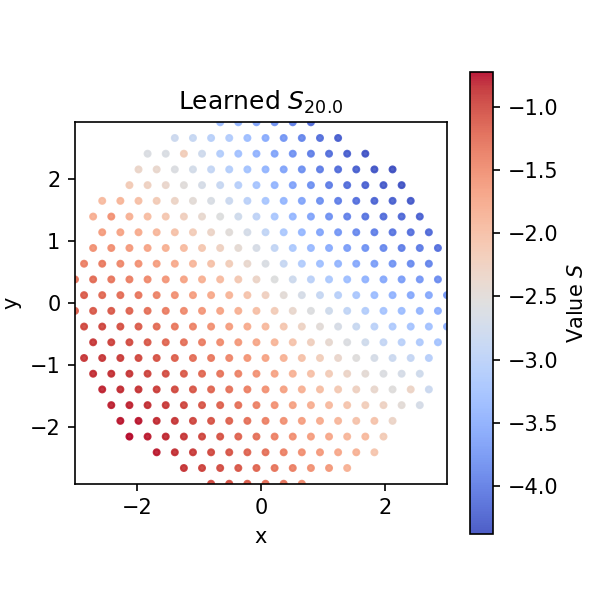}\subcaption{$t=20$}
    \end{minipage}
  }
  \end{minipage}
  \caption{Result of the dynamics of population density $\rho$ (upper panel) and value function $S$ (lower panel) with Gaussian interaction potential on the regular graph with triangle grids and circle boundary.}\label{fig:re-Gaussian}
\end{figure*}

\end{example}

\begin{example}{Lattice Grids and Fisher Information}

In this example, the graph $G$ is a lattice graph embedded in a circle $\mathcal{O}=x^2+y^2<3$. It has 489 nodes and 928 undirected edges. 
The initial density of the crowd is $\mu_0(x) = \rho_G (x,\mathbf{-0.75},0.2\cdot \mathbf{I}_2),$ and the target is $ \mu_T(x) = \rho_G (x,\mathbf{0.75},0.2\cdot \mathbf{I}_2).$ We consider Fisher information in this example, where $I(\rho_i)=\frac{\delta \mathcal{I}(\rho)}{\delta \rho_i} = \sum_{j \in N(i)}[(\log\rho_i-\log\rho_j)\frac{\theta_{ij}(\rho)}{\rho_i} + \frac{1}{2}(\log\rho_i-\log\rho_j)^2\frac{\partial \theta_{ij}(\rho)}{\partial \rho_i}]$ with the weight $\lambda_I$. Noted that the first variation of Fisher information is locally Lipschitz continuous, we assume the solution exists here. In the loss function, the coefficient ratio is $\lambda_K:\lambda_I:\lambda_G=0.5:0.0001:5000.$ 
We apply 3000 epochs in the warm up training, with the linear combination parameter $\alpha=0.6$ in~\cref{eq:warmup}. 
The experimental settings are $T=25$, $\Delta t=0.5$ and $M=50$. The training settings are $lr=10^{-5}$ and training epochs to be 6000. 

\cref{fig:re-fisher} shows the dynamics of the population density and value function. The crowd spreads out a bit even until it reach to the target $\mu_T$. In the result, the total kinetic energy is 7.7311, Fisher information is 0.0760, terminal energy reach 0.5339, and $\frac{1}{n}\|\rho(T)-\mu_T\|_1=0.0011$. We perform a comparison experiment that is without Fisher information in the loss. In this case, the value of the Fisher information of the predicted dynamics is 0.1087, with a smaller kinetic energy 5.9003, and smaller terminal energy 0.0326, which indicates that the population has a better match to the target density. 

\begin{figure*}[htbp] 
  \centering
  \begin{minipage}[b]{0.99\textwidth} 
  {
    \begin{minipage}[b]{0.15\linewidth} 
      \centering
      \includegraphics[width=\linewidth]{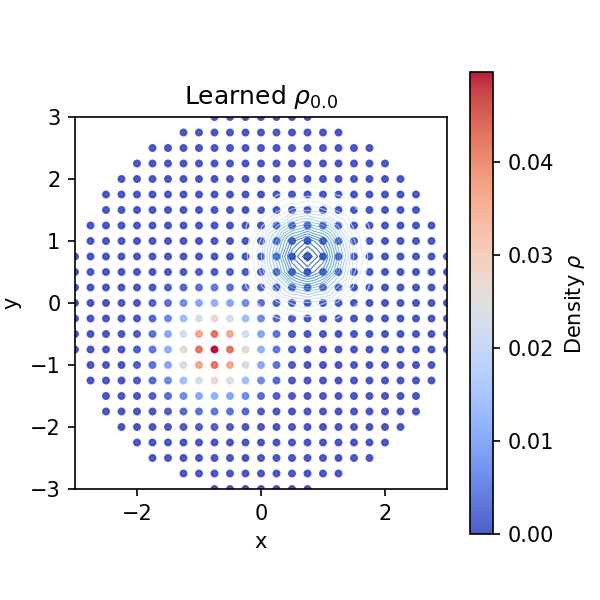}
      \includegraphics[width=\linewidth]{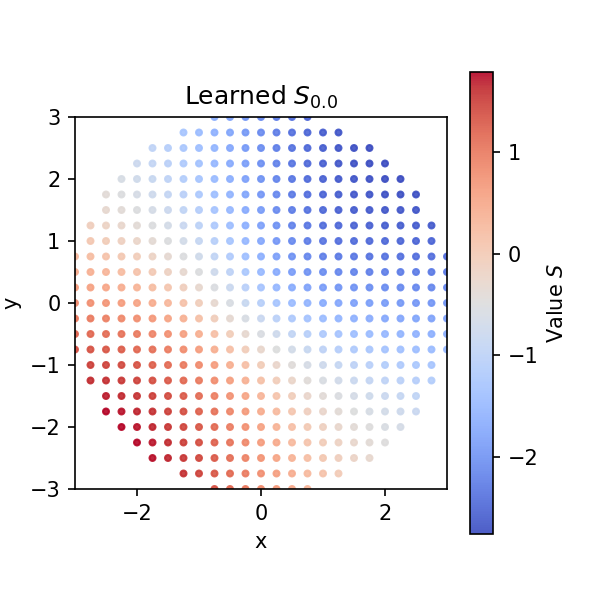}\subcaption{$t=0$}
    \end{minipage}
  }
   {
    \begin{minipage}[b]{0.15\linewidth}
      \centering
      \includegraphics[width=\linewidth]{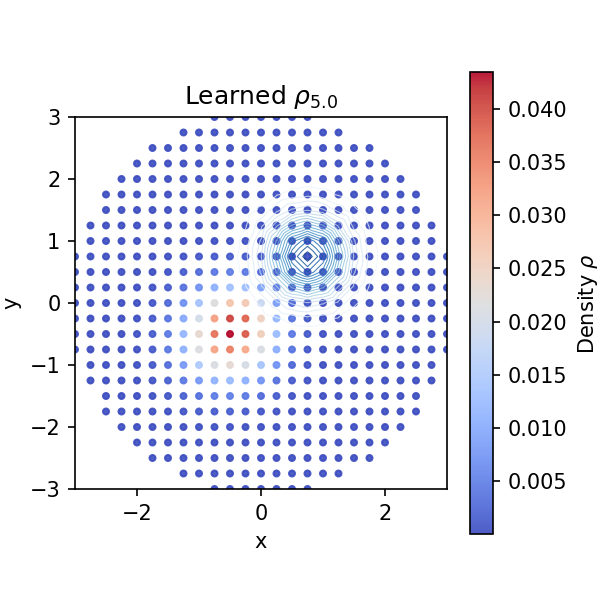}
      \includegraphics[width=\linewidth]{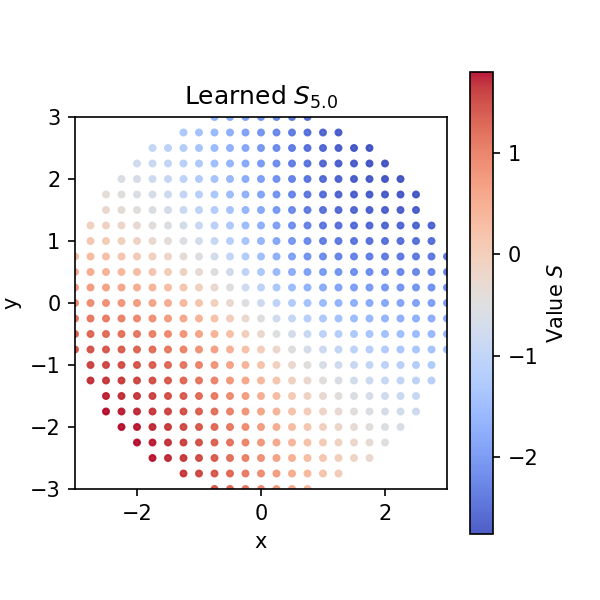}\subcaption{$t=5$}
    \end{minipage}
  }
   {
    \begin{minipage}[b]{0.15\linewidth}
      \centering
      \includegraphics[width=\linewidth]{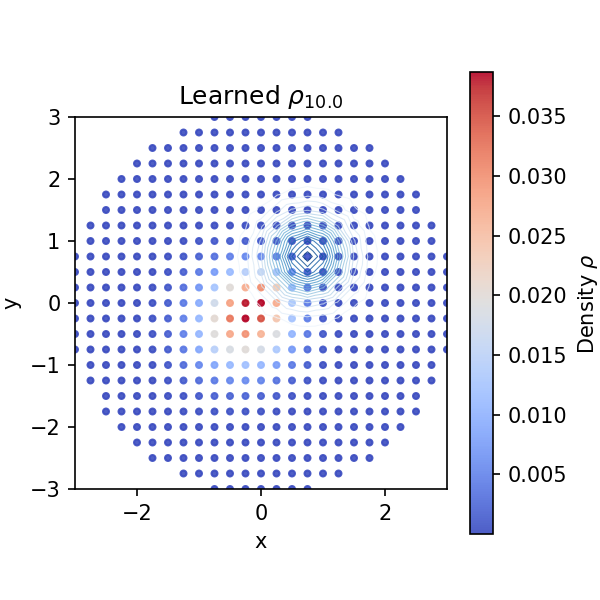}
      \includegraphics[width=\linewidth]{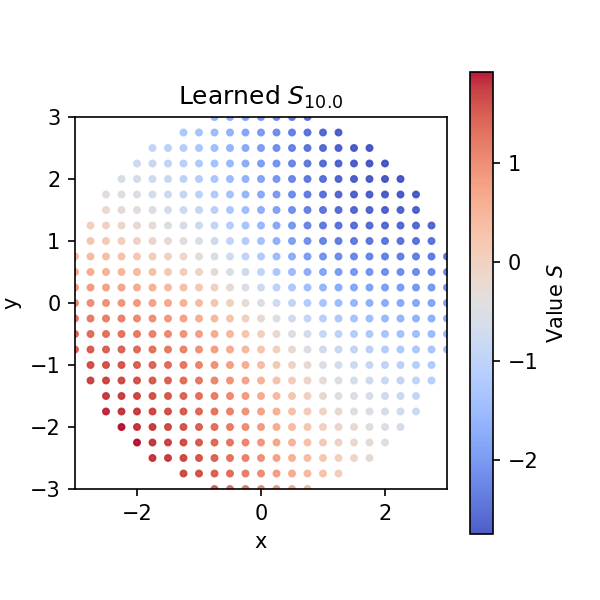}\subcaption{$t=10$}
    \end{minipage}
  }
    {
    \begin{minipage}[b]{0.15\linewidth}
      \centering
      \includegraphics[width=\linewidth]{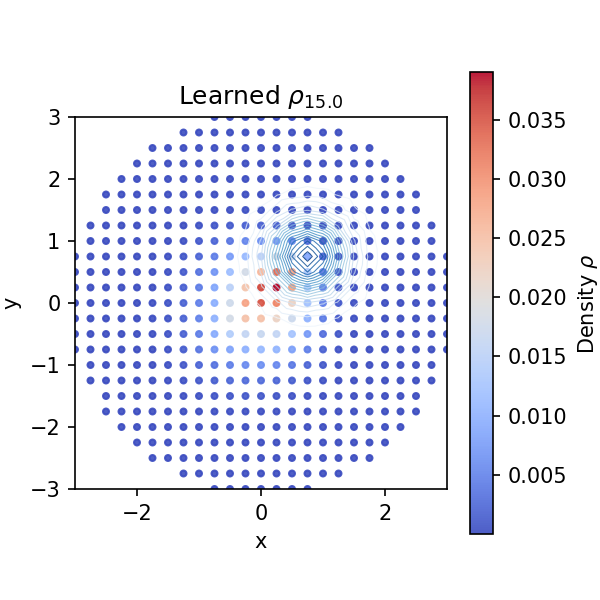}
      \includegraphics[width=\linewidth]{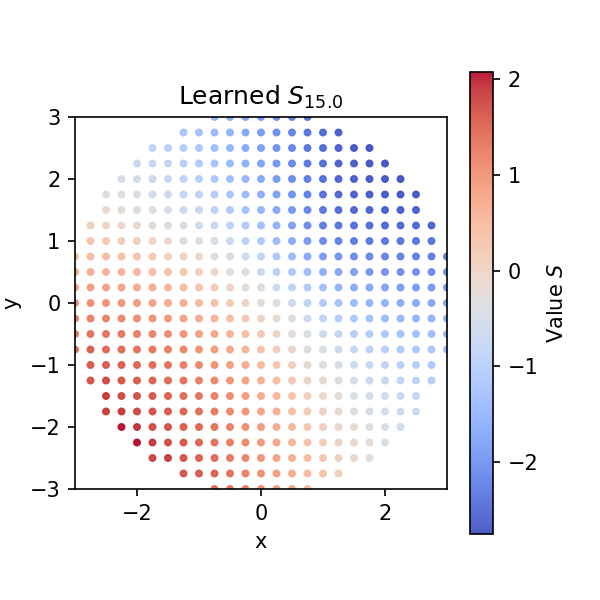}\subcaption{$t=15$}
    \end{minipage}
  }
    {
    \begin{minipage}[b]{0.15\linewidth}
      \centering
      \includegraphics[width=\linewidth]{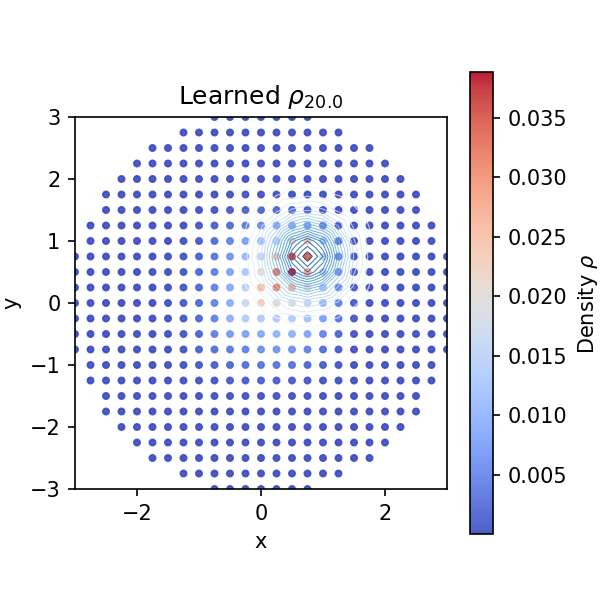}
      \includegraphics[width=\linewidth]{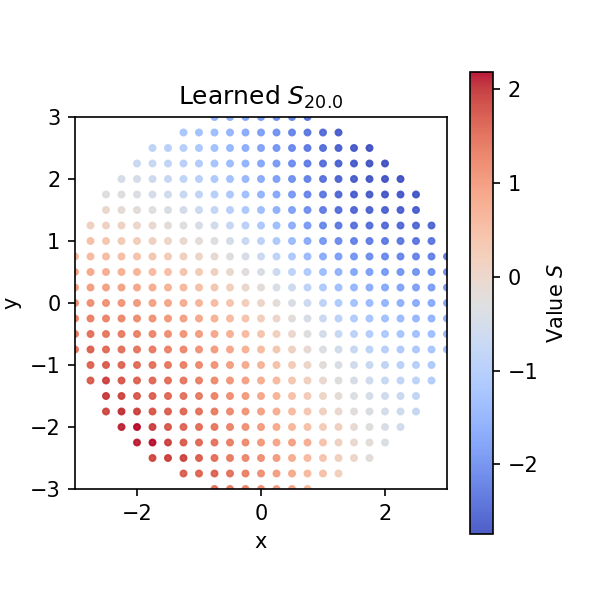}\subcaption{$t=20$}
    \end{minipage}
    \begin{minipage}[b]{0.15\linewidth}
   	\centering
   	\includegraphics[width=\linewidth]{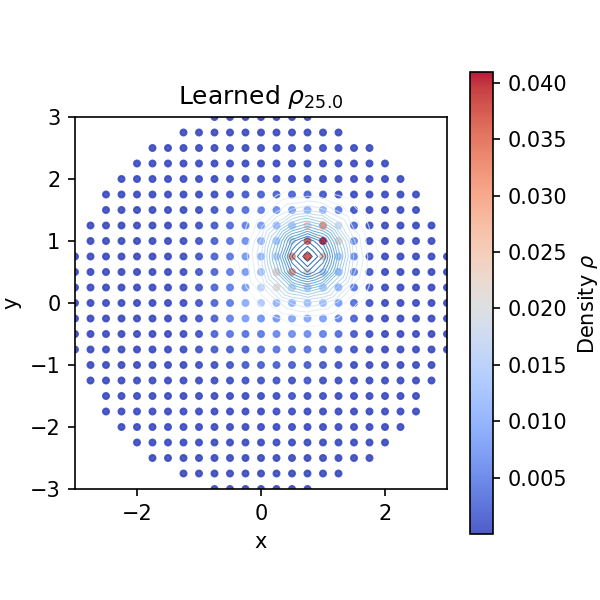}
   	\includegraphics[width=\linewidth]{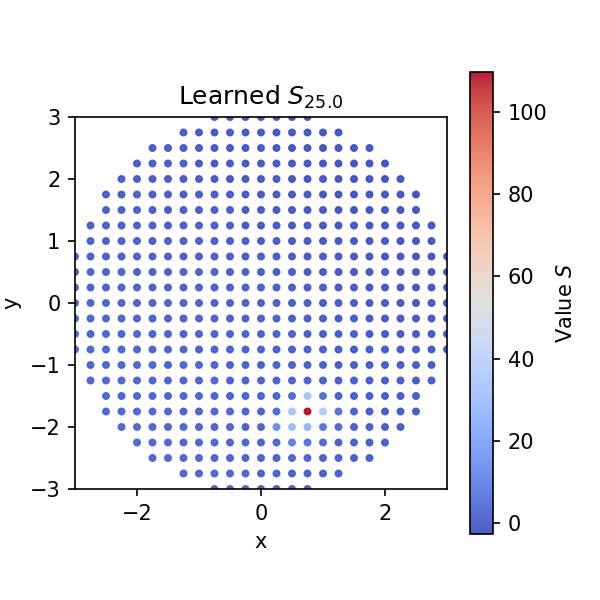}\subcaption{$t=25$}
    \end{minipage}
  }
  \end{minipage}
  \caption{Result of the dynamics of population density $\rho$ (upper panel) and value function $S$ (lower panel) with Fisher information on the regular graph with lattice grids and circle boundary.}\label{fig:re-fisher}
\end{figure*}

\end{example}


\begin{example}{Lattice Grids and Negative Boltzmann-Shannon Entropy}

In this example, the graph $G$ is a lattice graph embedded in a circle $\mathcal{O}=x^2+y^2<3$. It has 697 nodes and 1336 undirected edges.
The initial density of the crowd is given as $\mu_0(x) = \rho_G (x,\mathbf{-0.75},0.2\cdot \mathbf{I}_2),$ and the target is $ \mu_T(x) = \rho_G (x,\mathbf{0.75},0.2\cdot \mathbf{I}_2).$ Noted that the first variation of Boltzmann-Shannon entropy is locally Lipschitz continuous, we assume the solution exists here. In the loss function, the coefficient ratio is $\lambda_K:\lambda_B:\lambda_G=0.5:0.0001:5000.$ 
The experimental settings are $T=50$, $\Delta t=0.5$ and $M=100$. 

In the result, the entropy is -884.1543, terminal energy reaches 0.0689, and $\frac{1}{n}\|\rho(T)-\mu_T\|_1=9.8907 \times 10^{-5}$. In order to compare the effect brought by the entropy, we perform an experiment without the entropy in the loss, in which the value of the computed entropy converges to -872.8578, with a similar kinetic energy 4.3569, and smaller terminal energy 0.0144 that better arrival to the target density. 

\end{example}

\subsubsection{Graph MFG on Inhomogeneous Graphs}

\begin{example}{Inhomogeneous Graphs and Columbia Interaction Potential}

In this example, the graph $G$ is an inhomogeneous graph embedded in a circle $\mathcal{O}=x^2+y^2<3$. It has 2000 nodes and 8002 undirected edges.
The nodes on this randomly generated graph has the number of neighborhoods ranging from 6 to 10. 
The initial density of the crowd is $\mu_0(x) = \rho_G (x,-1.2 \mathbf{e}_1 + 1.3\mathbf{e}_2,0.2\cdot \mathbf{I}_2),$ and the target is $\mu_T(x) = \text{exp}(-a_0|x_0-b_0|-a_1|x_1-b_1|),$ in which $a_0=a_1=2, b_0=1, b_1=-1.25.$ In the loss function, the coefficient ratio is $\lambda_K:\lambda_W:\lambda_G=0.5:0.1:5000.$ The Columbia interaction potential is $\mathcal{W}(\rho):=\frac{1}{2} \sum_{j=1}^n \sum_{i=1}^n \frac{1}{\|X_i-Y_j\|^2+c} \rho_i \rho_j,$
where $c=0.5$. 
The experimental settings are $T=40$, $\Delta t=0.8$ and $M=50$.

\Cref{fig:re-columbia-irregular} shows the dynamics of the population density and value function. In this case, the inhomogeneous graph is randomly generated, unlike the regular graphs with fix neighbor nodes in the previous examples. In the dynamics, the crowd spreads out from the center and successfully arrives at the target. The spread out is due to the  Columbia interaction potential that penalizes aggregation. 
In the result, the total kinetic energy is 1.3480, interaction cost is 0.1826, terminal energy is 0.5312, and $\frac{1}{n}\|\rho(T)-\mu_T\|_1=0.0003$. 

\begin{figure*}[htbp] 
	\centering
	\begin{minipage}[b]{0.99\textwidth} 
		{
		\begin{minipage}[b]{0.14\linewidth} 
		\centering
		\includegraphics[width=\linewidth]{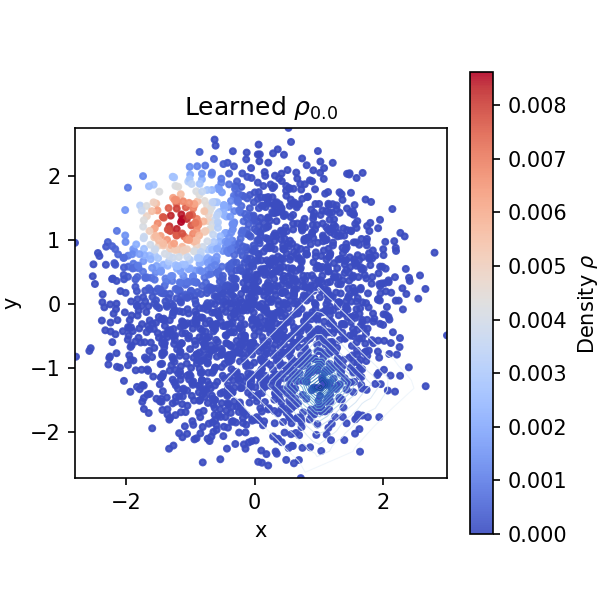}
		\includegraphics[width=\linewidth]{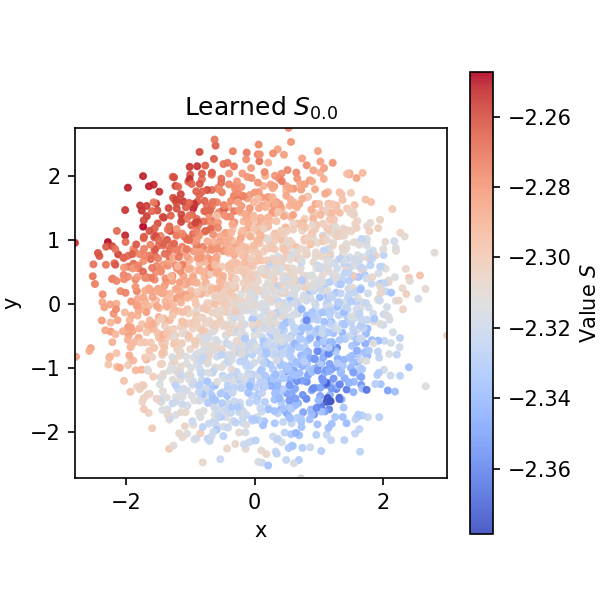}\subcaption{$t=0$}
		\end{minipage}
		}
		{
		\begin{minipage}[b]{0.14\linewidth}
		\centering
		\includegraphics[width=\linewidth]{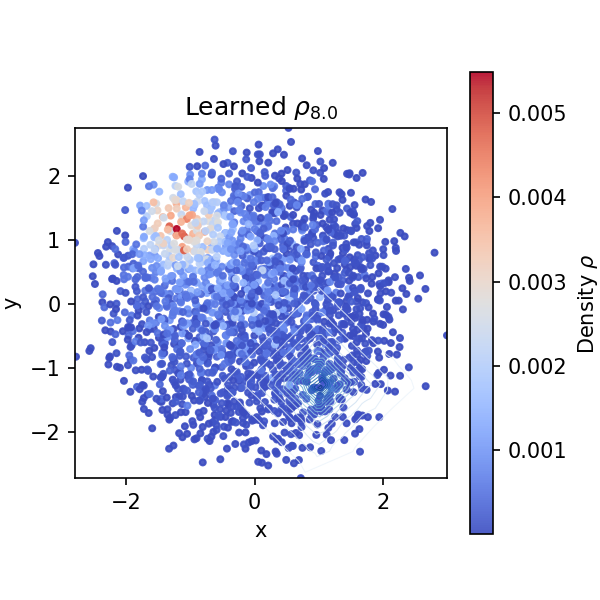}
		\includegraphics[width=\linewidth]{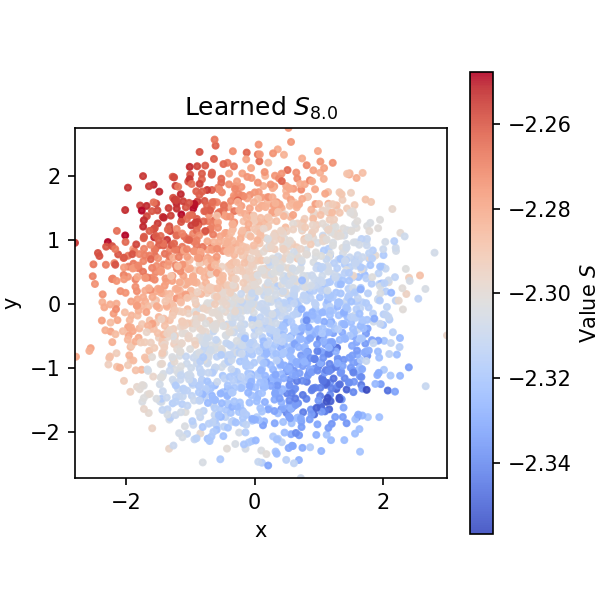}\subcaption{$t=8$}
		\end{minipage}
		}
		{
		\begin{minipage}[b]{0.14\linewidth}
		\centering
		\includegraphics[width=\linewidth]{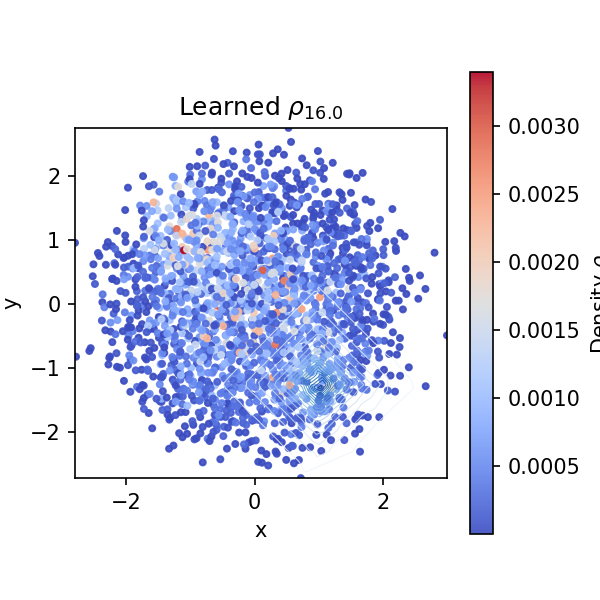}
		\includegraphics[width=\linewidth]{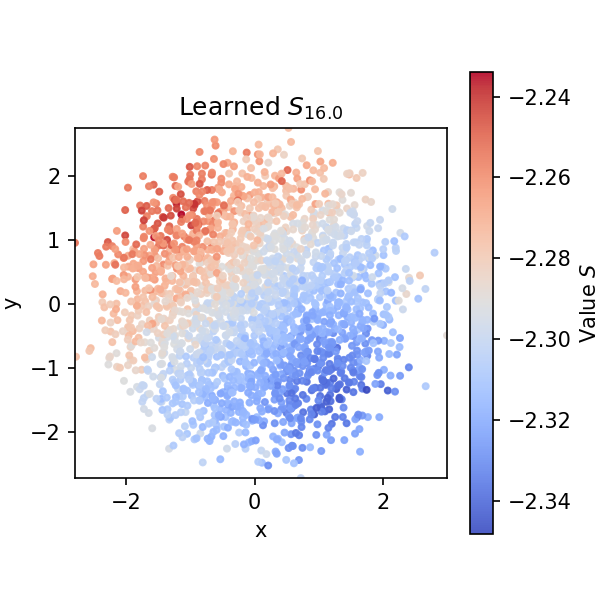}\subcaption{$t=16$}
		\end{minipage}
		}
		{
		\begin{minipage}[b]{0.14\linewidth}
		\centering
		\includegraphics[width=\linewidth]{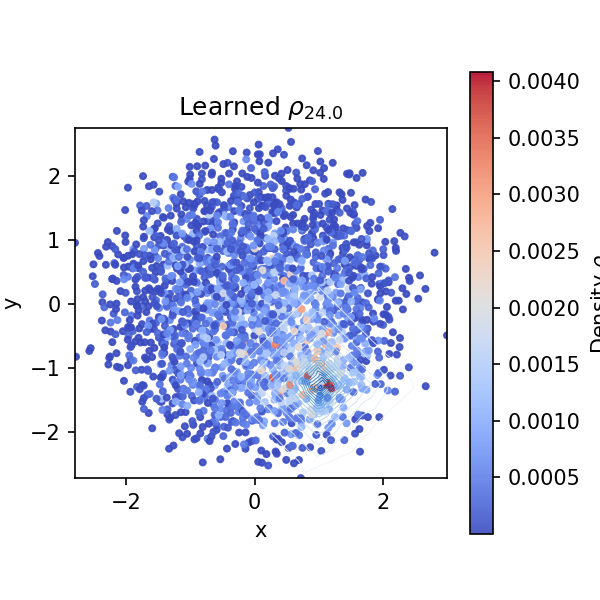}
		\includegraphics[width=\linewidth]{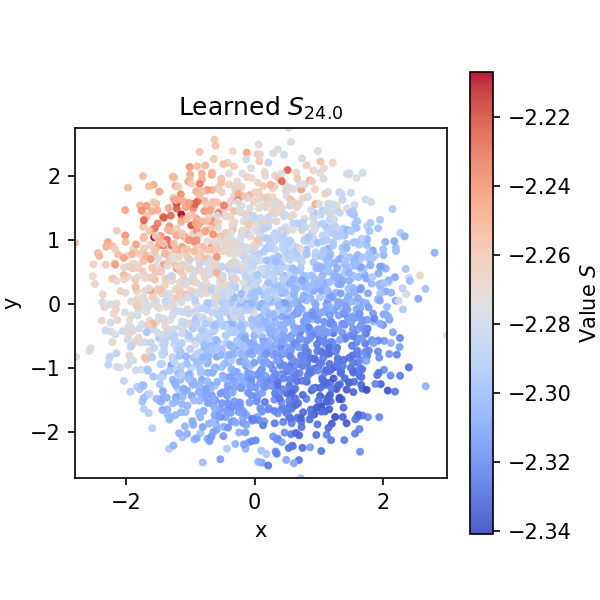}\subcaption{$t=24$}
		\end{minipage}
		}
		{
		\begin{minipage}[b]{0.14\linewidth}
		\centering
		\includegraphics[width=\linewidth]{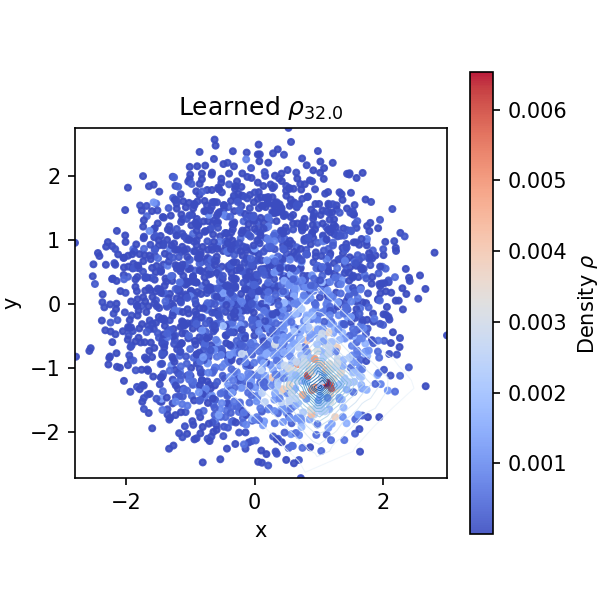}
		\includegraphics[width=\linewidth]{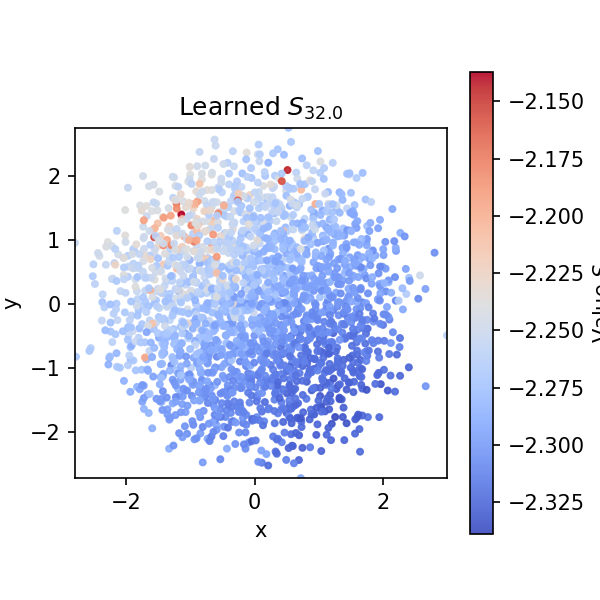}\subcaption{$t=32$}
		\end{minipage}
		}
		{
		\begin{minipage}[b]{0.14\linewidth}
		\centering
		\includegraphics[width=\linewidth]{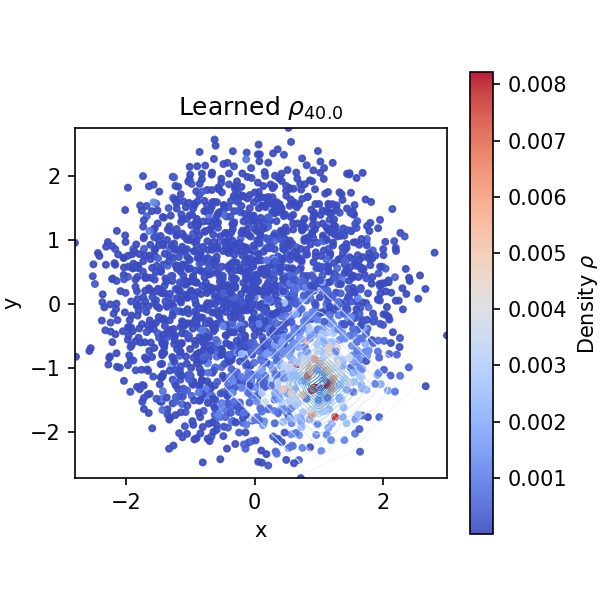}
		\includegraphics[width=\linewidth]{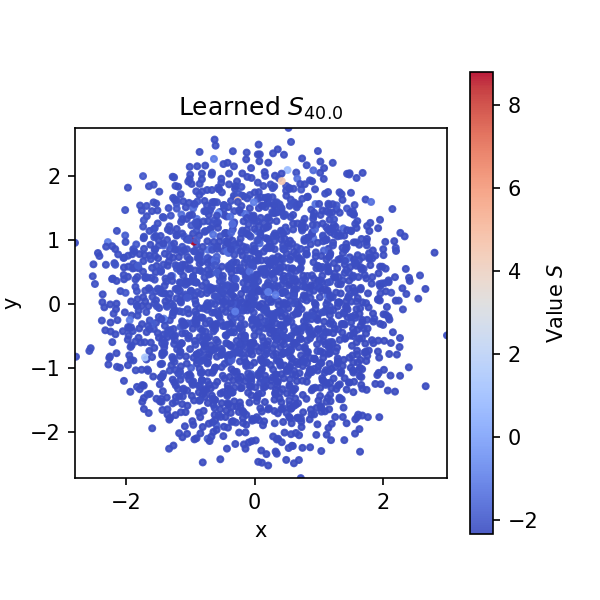}\subcaption{$t=40$}
		\end{minipage}
		}
	\end{minipage}
  \caption{Result of the dynamics of population dynamics $\rho$ (upper panel) and value function $S$ (lower panel) with Columbia interaction potential on the inhomogeneous graph, on which every node has the number of neighbor nodes ranging from 6 to 10.}\label{fig:re-columbia-irregular}
\end{figure*}

\end{example}

\begin{example}{Inhomogeneous Graph and Bifurcation}

In this example, the graph $G$ is an inhomogeneous graph embedded in a circle $\mathcal{O}=x^2+y^2<3$. It has 1000 nodes and 2747 undirected edges.
The nodes on this randomly generated graph has the number of neighborhoods ranging from 3 to 8. 
The initial distribution of the crowd is a Laplacian distribution, $\mu_0(x) = \text{exp}(-a_0|x_0-b_0|-a_1|x_1-b_1|),$ in which $a_0=a_1=2, b_0=b_1=0.$ The target distribution is a Gaussian mixture, that $\mu_T (x) =\frac{1}{6} \sum_{i=1}^6 \rho_G(x,\mathbf{m}_i,0.1\cdot \mathbf{I}_2),$ where $\mathbf{m}_i = \mathbf{p}_0 + r(sin(\theta_i)\mathbf{e}_1+cos(\theta_n))\mathbf{e}_2, \mathbf{p}_0=\mathbf{0}, r=1.5, \theta_i=\frac{2\pi}{6} i, i=1,\cdots,6.$ 
The experimental settings are $T=25$, $\Delta t=0.25$ and $M=100$. 

\Cref{fig:re-bifurcation-irregular} shows the dynamics of the population density and value function. In this case, a bifurcation of the crowd appears from the center of the graph to six groups of Gaussians. 
In the result, the total kinetic energy is 0.1697, terminal energy is 0.2413, and $\frac{1}{n}\|\rho(T)-\mu_T\|_1=0.0007$.


\begin{figure*}[htbp] 
  \centering
  \begin{minipage}[b]{0.99\textwidth} 
  {
    \begin{minipage}[b]{0.14\linewidth} 
      \centering
      \includegraphics[width=\linewidth]{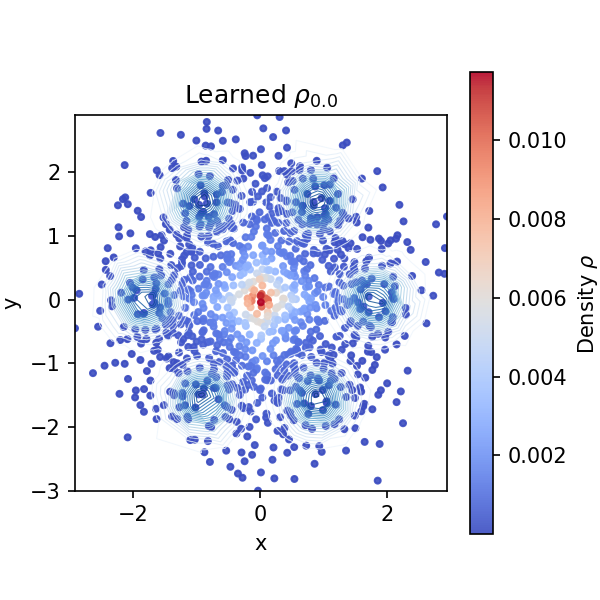}
      \includegraphics[width=\linewidth]{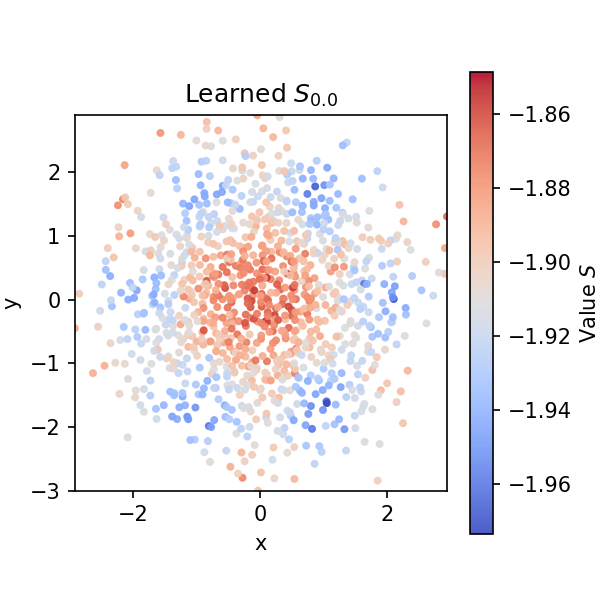}\subcaption{$t=0$}
    \end{minipage}
  }
   {
    \begin{minipage}[b]{0.14\linewidth}
      \centering
      \includegraphics[width=\linewidth]{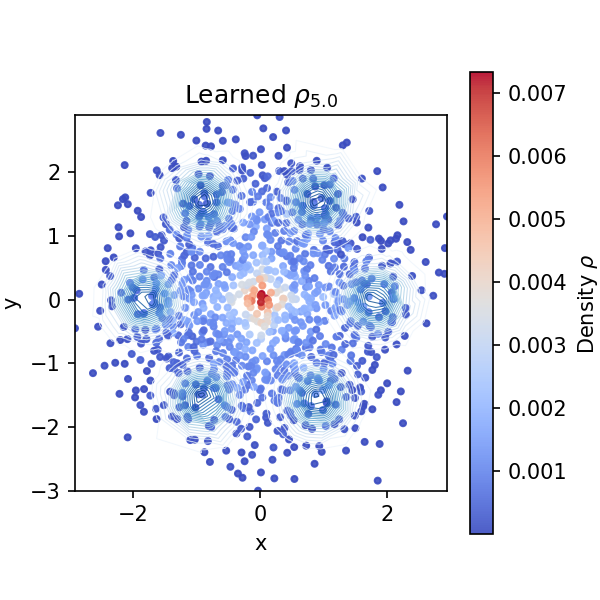}
      \includegraphics[width=\linewidth]{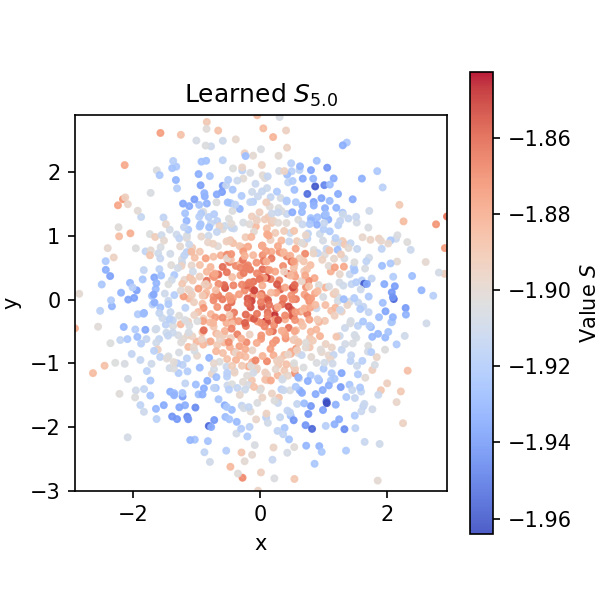}\subcaption{$t=5$}
    \end{minipage}
  }
   {
    \begin{minipage}[b]{0.14\linewidth}
      \centering
      \includegraphics[width=\linewidth]{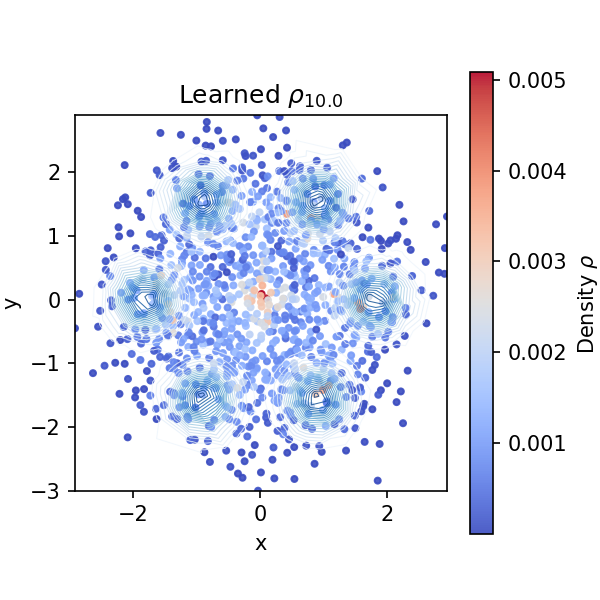}
      \includegraphics[width=\linewidth]{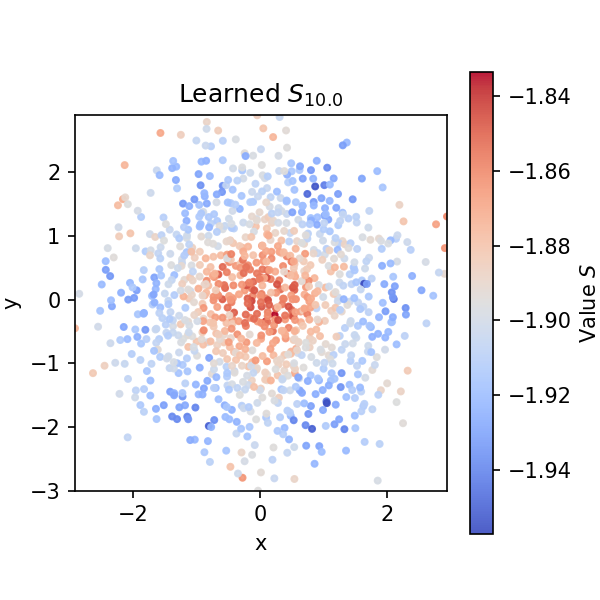}\subcaption{$t=10$}
    \end{minipage}
  }
    {
    \begin{minipage}[b]{0.14\linewidth}
      \centering
      \includegraphics[width=\linewidth]{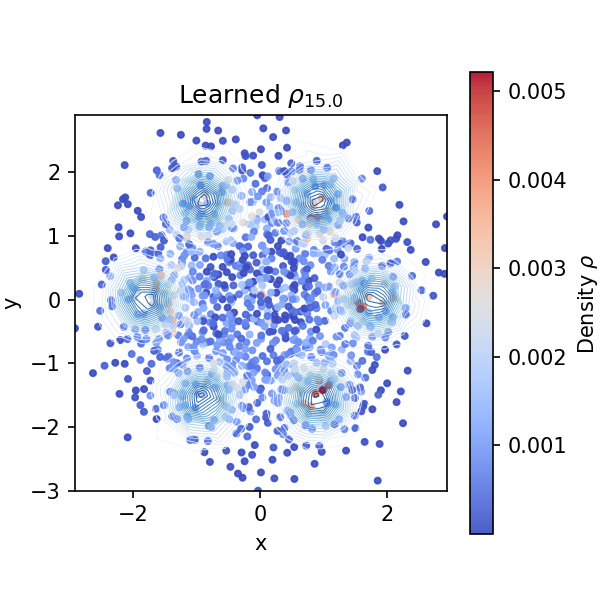}
      \includegraphics[width=\linewidth]{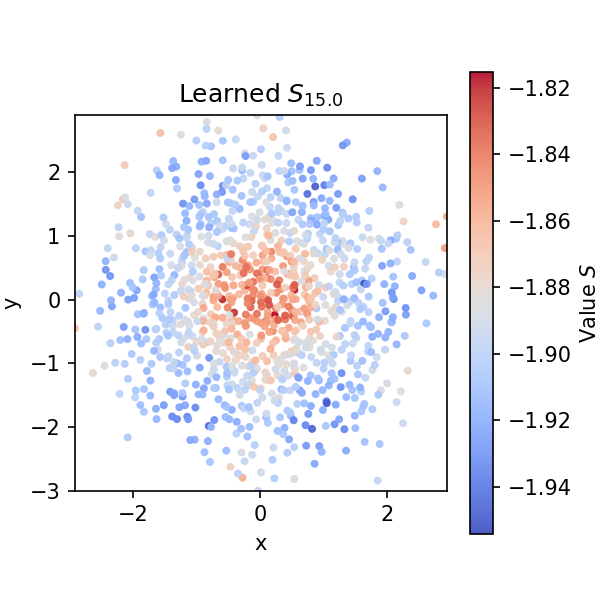}\subcaption{$t=15$}
    \end{minipage}
  }
    {
    \begin{minipage}[b]{0.14\linewidth}
      \centering
      \includegraphics[width=\linewidth]{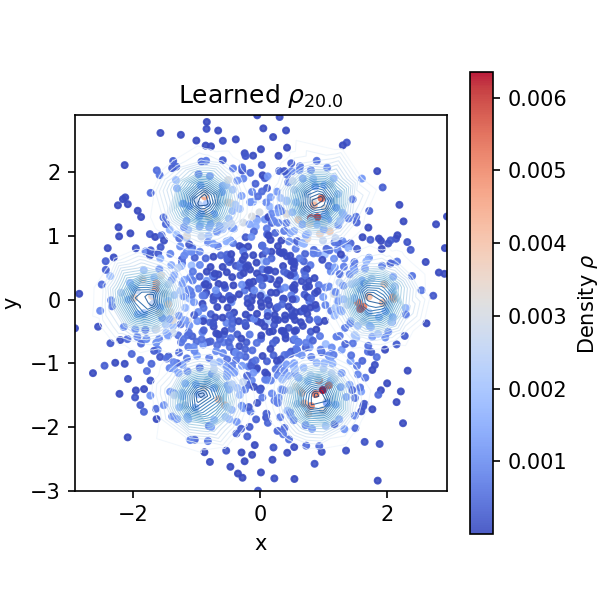}
      \includegraphics[width=\linewidth]{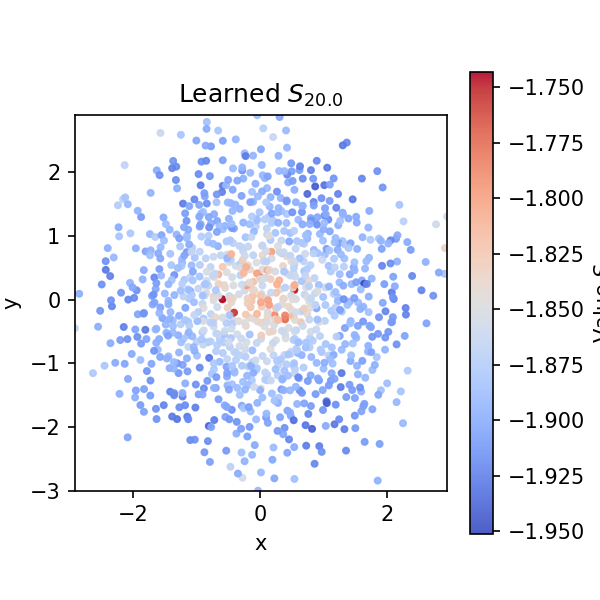}\subcaption{$t=20$}
    \end{minipage}
  }
      {
    \begin{minipage}[b]{0.14\linewidth}
      \centering
      \includegraphics[width=\linewidth]{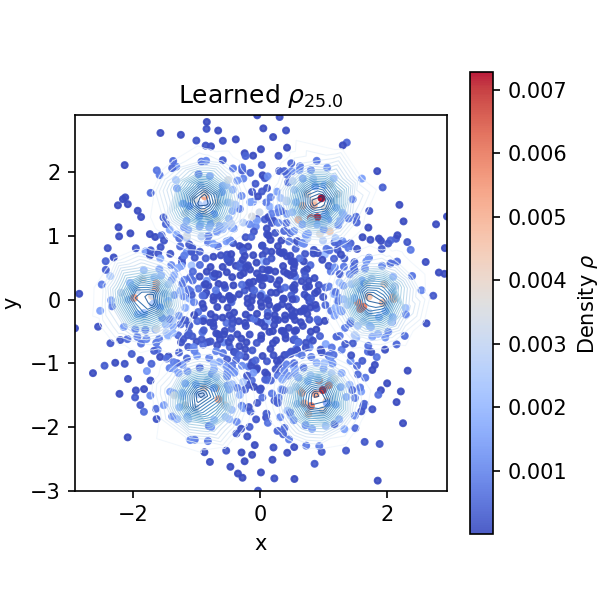}
      \includegraphics[width=\linewidth]{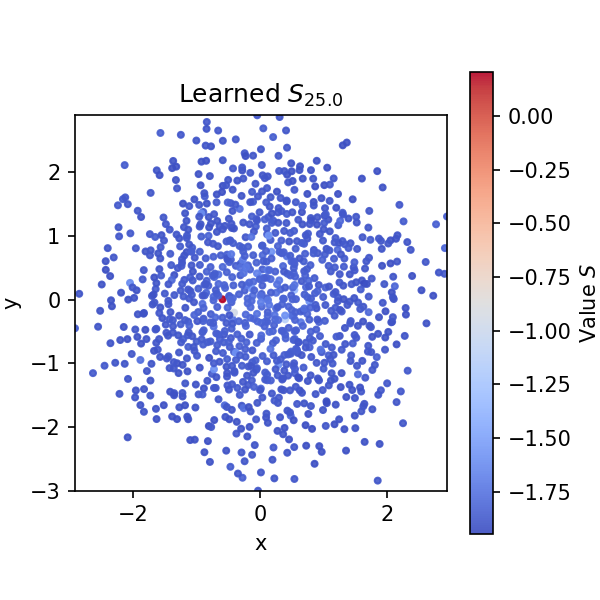}\subcaption{$t=25$}
    \end{minipage}
  }
  \end{minipage}
  \caption{Result of the dynamics of population dynamics $\rho$ (upper panel) and value function $S$ (lower panel) on the inhomogeneous graph, on which every node has the number of neighbor nodes ranging from 3 to 8.}\label{fig:re-bifurcation-irregular}
\end{figure*}

\end{example}

\section{Conclusions}
\label{sec:conclusions}
In this paper, we study the discrete potential mean field games on finite graphs. We propose an initial value optimization formulation Graph MFG-IV, which is finite-dimensional, and has a much smaller searching space than the path-wise formulation Graph MFG-PW. We prove the equivalent between these two formulations. Based on this formulation, we design a neural network algorithm to solve it. A warm-up scheme is proposed to improve the training. We apply our method to various Graph MFG examples with different graphs including regular grids, irregular boundaries, and inhomogeneous graphs, and different potentials including linear, nonlinear, and interaction ones. 
Because our approach is suitable for sequential processes as discrete mass transport models, it could be applied to the generative AI and Large Language Model in the future, such as the generation of text and code.



\section*{Acknowledgments}
This work was done during YXF's internship under the supervision of Professor Zhou, in the School of Mathematics, Georgia Institute of Technology. YXF and YX was partially supported by the Project of Hetao Shenzhen-HKUST Innovation Cooperation Zone HZQB-KCZYB-2020083.

\bibliographystyle{siamplain}
\bibliography{references}

@article{ruthotto2020machine,
  title={A machine learning framework for solving high-dimensional mean field game and mean field control problems},
  author={Ruthotto, Lars and Osher, Stanley J and Li, Wuchen and Nurbekyan, Levon and Fung, Samy Wu},
  journal={Proceedings of the National Academy of Sciences},
  volume={117},
  number={17},
  pages={9183--9193},
  year={2020},
  publisher={National Academy of Sciences}
}

@article{lasry2007mean,
  title={Mean field games},
  author={Lasry, Jean-Michel and Lions, Pierre-Louis},
  journal={Japanese journal of mathematics},
  volume={2},
  number={1},
  pages={229--260},
  year={2007},
  publisher={Springer}
}

@book{gross2003handbook,
  title={Handbook of graph theory},
  author={Gross, Jonathan L and Yellen, Jay},
  year={2003},
  publisher={CRC press}
}

@article{chow2022dynamical,
  title={Dynamical Schr{\"o}dinger bridge problems on graphs},
  author={Chow, Shui-Nee and Li, Wuchen and Mou, Chenchen and Zhou, Haomin},
  journal={Journal of Dynamics and Differential Equations},
  volume={34},
  number={3},
  pages={2511--2530},
  year={2022},
  publisher={Springer}
}

@article{cui2022time,
  title={Time discretizations of Wasserstein--Hamiltonian flows},
  author={Cui, Jianbo and Dieci, Luca and Zhou, Haomin},
  journal={Mathematics of Computation},
  volume={91},
  number={335},
  pages={1019--1075},
  year={2022}
}

@article{chow2019discrete,
  title={A discrete Schr{\"o}dinger equation via optimal transport on graphs},
  author={Chow, Shui-Nee and Li, Wuchen and Zhou, Haomin},
  journal={Journal of Functional Analysis},
  volume={276},
  number={8},
  pages={2440--2469},
  year={2019},
  publisher={Elsevier}
}

@article{chow2017entropy,
  title={Entropy dissipation of Fokker-Planck equations on graphs},
  author={Chow, Shui-Nee and Li, Wuchen and Zhou, Haomin},
  journal={arXiv preprint arXiv:1701.04841},
  year={2017}
}

@article{hamilton2017inductive,
  title={Inductive representation learning on large graphs},
  author={Hamilton, Will and Ying, Zhitao and Leskovec, Jure},
  journal={Advances in neural information processing systems},
  volume={30},
  year={2017}
}

@article{scarselli2008graph,
  title={The graph neural network model},
  author={Scarselli, Franco and Gori, Marco and Tsoi, Ah Chung and Hagenbuchner, Markus and Monfardini, Gabriele},
  journal={IEEE transactions on neural networks},
  volume={20},
  number={1},
  pages={61--80},
  year={2008},
  publisher={IEEE}
}

@article{hu2023recent,
  title={Recent developments in machine learning methods for stochastic control and games},
  author={Hu, Ruimeng and Lauriere, Mathieu},
  journal={arXiv preprint arXiv:2303.10257},
  year={2023}
}

@article{liu2022deep,
  title={Deep generalized schr{\"o}dinger bridge},
  author={Liu, Guan-Horng and Chen, Tianrong and So, Oswin and Theodorou, Evangelos},
  journal={Advances in Neural Information Processing Systems},
  volume={35},
  pages={9374--9388},
  year={2022}
}

@article{lin2021alternating,
  title={Alternating the population and control neural networks to solve high-dimensional stochastic mean-field games},
  author={Lin, Alex Tong and Fung, Samy Wu and Li, Wuchen and Nurbekyan, Levon and Osher, Stanley J},
  journal={Proceedings of the National Academy of Sciences},
  volume={118},
  number={31},
  pages={e2024713118},
  year={2021},
  publisher={National Academy of Sciences}
}

@article{lawal2022physics,
  title={Physics-informed neural network (PINN) evolution and beyond: A systematic literature review and bibliometric analysis},
  author={Lawal, Zaharaddeen Karami and Yassin, Hayati and Lai, Daphne Teck Ching and Che Idris, Azam},
  journal={Big Data and Cognitive Computing},
  volume={6},
  number={4},
  pages={140},
  year={2022},
  publisher={MDPI}
}

@article{huang2023bridging,
  title={Bridging mean-field games and normalizing flows with trajectory regularization},
  author={Huang, Han and Yu, Jiajia and Chen, Jie and Lai, Rongjie},
  journal={Journal of Computational Physics},
  volume={487},
  pages={112155},
  year={2023},
  publisher={Elsevier}
}

@article{gangbo2024well,
  title={Well-posedness for Hamilton--Jacobi equations on the Wasserstein space on graphs},
  author={Gangbo, Wilfrid and Mou, Chenchen and {\'S}wi{\k{e}}ch, Andrzej},
  journal={Calculus of Variations and Partial Differential Equations},
  volume={63},
  number={7},
  pages={160},
  year={2024},
  publisher={Springer}
}

@article{cui2025finite,
  title={Finite difference schemes for Hamilton--Jacobi equation on Wasserstein space on graphs},
  author={Cui, Jianbo and Dang, Tonghe and Mou, Chenchen},
  journal={arXiv preprint arXiv:2504.13463},
  year={2025}
}

@article{han2024learning,
  title={Learning high-dimensional McKean--Vlasov forward-backward stochastic differential equations with general distribution dependence},
  author={Han, Jiequn and Hu, Ruimeng and Long, Jihao},
  journal={SIAM Journal on Numerical Analysis},
  volume={62},
  number={1},
  pages={1--24},
  year={2024},
  publisher={SIAM}
}

@article{benamou2000computational,
  title={A computational fluid mechanics solution to the Monge-Kantorovich mass transfer problem},
  author={Benamou, Jean-David and Brenier, Yann},
  journal={Numerische Mathematik},
  volume={84},
  number={3},
  pages={375--393},
  year={2000},
  publisher={Springer-Verlag Berlin/Heidelberg}
}

@article{carmona2022convergence,
  title={Convergence analysis of machine learning algorithms for the numerical solution of mean field control and games: II—the finite horizon case},
  author={Carmona, Ren{\'e} and Lauri{\`e}re, Mathieu},
  journal={The Annals of Applied Probability},
  volume={32},
  number={6},
  pages={4065--4105},
  year={2022},
  publisher={Institute of Mathematical Statistics}
}

@article{chow2012fokker,
  title={Fokker--Planck equations for a free energy functional or Markov process on a graph},
  author={Chow, Shui-Nee and Huang, Wen and Li, Yao and Zhou, Haomin},
  journal={Archive for Rational Mechanics and Analysis},
  volume={203},
  number={3},
  pages={969--1008},
  year={2012},
  publisher={Springer}
}

@article{chow2020wasserstein,
  title={Wasserstein hamiltonian flows},
  author={Chow, Shui-Nee and Li, Wuchen and Zhou, Haomin},
  journal={Journal of Differential Equations},
  volume={268},
  number={3},
  pages={1205--1219},
  year={2020},
  publisher={Elsevier}
}

@article{achdou2022income,
  title={Income and wealth distribution in macroeconomics: A continuous-time approach},
  author={Achdou, Yves and Han, Jiequn and Lasry, Jean-Michel and Lions, Pierre-Louis and Moll, Benjamin},
  journal={The review of economic studies},
  volume={89},
  number={1},
  pages={45--86},
  year={2022},
  publisher={Oxford University Press}
}

@article{gueant2010mean,
  title={Mean field games and applications. paris-princeton lectures on mathematical finance},
  author={Gu{\'e}ant, O and Lasry, JM and Lions, PL},
  journal={Lect. Notes Math},
  volume={2011},
  pages={205--266},
  year={2010}
}

@article{carmona2020applications,
  title={Applications of mean field games in financial engineering and economic theory},
  author={Carmona, Rene},
  journal={arXiv preprint arXiv:2012.05237},
  year={2020}
}

@article{bauso2016opinion,
  title={Opinion dynamics in social networks through mean-field games},
  author={Bauso, Dario and Tembine, Hamidou and Basar, Tamer},
  journal={SIAM Journal on Control and Optimization},
  volume={54},
  number={6},
  pages={3225--3257},
  year={2016},
  publisher={SIAM}
}

@incollection{cardaliaguet2021introduction,
  title={An introduction to mean field game theory},
  author={Cardaliaguet, Pierre and Porretta, Alessio},
  booktitle={Mean Field Games: Cetraro, Italy 2019},
  pages={1--158},
  year={2021},
  publisher={Springer}
}

@book{chew2016potential,
  title={Potential game theory},
  author={Chew, Yong Huat and Soong, Boon-Hee and others},
  year={2016},
  publisher={Springer}
}

@inproceedings{liu2018mean,
  title={A mean field game approach to swarming robots control},
  author={Liu, Zhiyu and Wu, Bo and Lin, Hai},
  booktitle={2018 Annual American Control Conference (ACC)},
  pages={4293--4298},
  year={2018},
  organization={IEEE}
}

@article{elamvazhuthi2019mean,
  title={Mean-field models in swarm robotics: A survey},
  author={Elamvazhuthi, Karthik and Berman, Spring},
  journal={Bioinspiration \& Biomimetics},
  volume={15},
  number={1},
  pages={015001},
  year={2019},
  publisher={IOP Publishing}
}

@inproceedings{kang2020joint,
  title={Joint task assignment and trajectory optimization for a mobile robot swarm by mean-field game},
  author={Kang, Yuhan and Liu, Siting and Lee, Wonjun and Zhang, Hongliang and Li, Wuchen and Han, Zhu},
  booktitle={Globecom 2020-2020 IEEE global communications conference},
  pages={1--6},
  year={2020},
  organization={IEEE}
}

@article{aurell2022optimal,
  title={Optimal incentives to mitigate epidemics: a Stackelberg mean field game approach},
  author={Aurell, Alexander and Carmona, Rene and Dayanikli, Gokce and Lauriere, Mathieu},
  journal={SIAM Journal on Control and Optimization},
  volume={60},
  number={2},
  pages={S294--S322},
  year={2022},
  publisher={SIAM}
}

@article{roy2023recent,
  title={Recent advances in modeling and control of epidemics using a mean field approach},
  author={Roy, Amal and Singh, Chandramani and Narahari, Y},
  journal={S{\=a}dhan{\=a}},
  volume={48},
  number={4},
  pages={207},
  year={2023},
  publisher={Springer}
}

@article{yang2017learning,
  title={Learning deep mean field games for modeling large population behavior},
  author={Yang, Jiachen and Ye, Xiaojing and Trivedi, Rakshit and Xu, Huan and Zha, Hongyuan},
  journal={arXiv preprint arXiv:1711.03156},
  year={2017}
}

@inproceedings{chen2023learning,
  title={Learning dual mean field games on graphs},
  author={Chen, Xu and Liu, Shuo and Di, Xuan},
  booktitle={the European Conference on Artificial Intelligence (ECAI 2023)},
  year={2023},
  organization={European Conference on Artificial Intelligence}
}

@article{huang2006large,
  title={Large population stochastic dynamic games: closed-loop McKean-Vlasov systems and the Nash certainty equivalence principle},
  author={Huang, Minyi and Malham{\'e}, Roland P and Caines, Peter E},
  journal={Communications in Information and Systems},
  year={2006}
}

@inproceedings{calderone2017markov,
  title={Markov decision process routing games},
  author={Calderone, Dan and Sastry, S Shankar},
  booktitle={Proceedings of the 8th International Conference on Cyber-Physical Systems},
  pages={273--279},
  year={2017}
}

@article{huang2021dynamic,
  title={Dynamic driving and routing games for autonomous vehicles on networks: A mean field game approach},
  author={Huang, Kuang and Chen, Xu and Di, Xuan and Du, Qiang},
  journal={Transportation Research Part C: Emerging Technologies},
  volume={128},
  pages={103189},
  year={2021},
  publisher={Elsevier}
}

@article{cabannes2021solving,
  title={Solving N-player dynamic routing games with congestion: a mean field approach},
  author={Cabannes, Theophile and Lauriere, Mathieu and Perolat, Julien and Marinier, Raphael and Girgin, Sertan and Perrin, Sarah and Pietquin, Olivier and Bayen, Alexandre M and Goubault, Eric and Elie, Romuald},
  journal={arXiv preprint arXiv:2110.11943},
  year={2021}
}

@article{shou2022multi,
  title={Multi-agent reinforcement learning for Markov routing games: A new modeling paradigm for dynamic traffic assignment},
  author={Shou, Zhenyu and Chen, Xu and Fu, Yongjie and Di, Xuan},
  journal={Transportation Research Part C: Emerging Technologies},
  volume={137},
  pages={103560},
  year={2022},
  publisher={Elsevier}
}

@article{tanaka2020linearly,
  title={Linearly solvable mean-field traffic routing games},
  author={Tanaka, Takashi and Nekouei, Ehsan and Pedram, Ali Reza and Johansson, Karl Henrik},
  journal={IEEE Transactions on Automatic Control},
  volume={66},
  number={2},
  pages={880--887},
  year={2020},
  publisher={IEEE}
}

@article{maas2011gradient,
  title={Gradient flows of the entropy for finite Markov chains},
  author={Maas, Jan},
  journal={Journal of Functional Analysis},
  volume={261},
  number={8},
  pages={2250--2292},
  year={2011},
  publisher={Elsevier}
}

@article{mielke2011gradient,
  title={A gradient structure for reaction--diffusion systems and for energy-drift-diffusion systems},
  author={Mielke, Alexander},
  journal={Nonlinearity},
  volume={24},
  number={4},
  pages={1329},
  year={2011},
  publisher={IOP Publishing}
}

@article{gomes2010discrete,
  title={Discrete time, finite state space mean field games},
  author={Gomes, Diogo A and Mohr, Joana and Souza, Rafael Rigao},
  journal={Journal de math{\'e}matiques pures et appliqu{\'e}es},
  volume={93},
  number={3},
  pages={308--328},
  year={2010},
  publisher={Elsevier}
}

@article{gomes2013continuous,
  title={Continuous time finite state mean field games},
  author={Gomes, Diogo A and Mohr, Joana and Souza, Rafael Rigao},
  journal={Applied Mathematics \& Optimization},
  volume={68},
  number={1},
  pages={99--143},
  year={2013},
  publisher={Springer}
}

@article{gueant2011infinity,
  title={From infinity to one: The reduction of some mean field games to a global control problem},
  author={Gu{\'e}ant, Olivier},
  journal={arXiv preprint arXiv:1110.3441},
  year={2011}
}

@article{gueant2015existence,
  title={Existence and uniqueness result for mean field games with congestion effect on graphs},
  author={Gu{\'e}ant, Olivier},
  journal={Applied Mathematics \& Optimization},
  volume={72},
  number={2},
  pages={291--303},
  year={2015},
  publisher={Springer}
}

@article{zhang2023mean,
  title={A mean-field games laboratory for generative modeling},
  author={Zhang, Benjamin J and Katsoulakis, Markos A},
  journal={arXiv preprint arXiv:2304.13534},
  year={2023}
}

\end{document}